\documentclass[12pt,reqno]{amsart}
\usepackage{amsfonts,amssymb,amsmath}

\DeclareMathOperator{\bmin}{\mathbf{min}}
\DeclareMathOperator{\bmax}{\mathbf{max}}

\newtheorem{thm}{Theorem}
\newtheorem{defn}[thm]{Definition}
\newtheorem{cor}[thm]{Corollary}

\newtheorem{prop}[thm]{Proposition}

\newtheorem{rem}[thm]{Remark}
\newtheorem{exmpl}[thm]{Example}

\newenvironment{sketch}{\vspace*{2eX}\noindent
{\em Sketch of proof.} }{\hfill$\square$ \\[2eX]}

\numberwithin{equation}{section} \numberwithin{thm}{section}

\title{Relative Blocking in Posets}

\date{}
\author{Andrey O. Matveev}
\thanks{Data-Center Co., RU-620034, Ekaterinburg, P.O.~Box~5, Russian~Federation}
\thanks{aomatveev@\{dc.ru, hotmail.com\}}

\begin{document}

\begin{abstract}
Poset-theoretic generalizations of set-theoretic committee
constructions are presented. The structure of the corresponding
subposets is described. Sequences of irreducible fractions
associated to the principal order ideals of finite bounded posets
are considered and those related to the Boolean lattices are
explored; it is shown that such sequences inherit all the familiar
properties of the Farey sequences.
\end{abstract}

\maketitle

\section{Introduction and preliminaries}

Various decision-making, recognition, and voting procedures rely, explicitly or implicitly, on the
cardinalities of finite sets and of their mutual intersections. Among mathematical constructions which
underlie those procedures are {\em blocking sets\/} ({\em covers}, {\em systems of representatives}, {\em
transversals})~(F\"{u}redi, 1988, and Chapter~8 of Gr\"{o}tschel~et~al., 1988), {\em committees}~(Khachai et
al., 2002), and {\em quorum systems\/} ({\em intersecting set systems}, {\em intersecting
hypergraphs})~(Colbourn et al., 2001, Loeb and Conway, 2000, and Naor and Wool, 1998); see also~Crama and
Hammer (in preparation).

The present paper is devoted to discussing questions concerning
mechanisms of blocking in finite posets that go back to
set-theoretic committees.

We refer the reader to~Chapter~3 of Stanley, 1997, for information
and terminology in the theory of posets.

Recall that a set $H$ is called a {\em blocking set\/} for a
nonempty family $\mathcal{G}=\{G_1,\ldots,G_m\}$ of nonempty
subsets of a finite set if it holds $|H\cap G_k|>0$, for each
$k\in\{1,\ldots,m\}$. The family of all inclusion-minimal blocking
sets for $\mathcal{G}$ is called the {\em blocker\/} of
$\mathcal{G}$, see, e.g., Chapter~8 of Gr\"{o}tschel et al., 1993.
Let $r$ be a rational number such that $0\leq r<1$. A set $H$ is
called an $r$-{\em committee\/} for $\mathcal{G}$ if it holds
$|H\cap G_k|>r\cdot|H|$, for each $k\in\{1,\ldots,m\}$, see, e.g.,
Khachai et al., 2002.

A family of subsets of a finite {\em ground set\/} is called a
{\em clutter\/} or a {\em Sperner family\/} if no set from that
family contains another. The empty clutter containing no subsets
of the ground set, and the clutter whose unique set is the empty
subset of the ground set, are called the {\em trivial clutters}.
The {\em blocker map\/} assigns to a nontrivial clutter its
blocker, and this map assigns to a trivial clutter the other
trivial clutter, see, e.g., Cordovil et al., 1991.

The set-theoretic blocker constructions are at the foundation of discrete mathematics, see, e.g.,
Cornu\'{e}jols, 2001, and Crama and Hammer (in preparation).

Since the clutters on a ground set are in one-to-one
correspondence with the antichains in the Boolean lattice of all
subsets of the ground set, the set-theoretic concepts of blocking
can be assigned poset-theoretic counterparts. The next natural
step consists in a passage from the Boolean lattices to arbitrary
finite bounded posets, see~Bj\"{o}rner et al., 2004, 2005, and
Matveev, 2001, 2002, 2003; a poset is called {\em bounded\/} if it
has a least and greatest elements.

Throughout the paper, $P$ stands for a finite bounded poset of
cardinality greater than one whose least and greatest elements are
denoted by $\hat{0}_P$ and $\hat{1}_P$, respectively. $P^{\mathrm
a}$ denotes the set of all atoms of $P$ (the {\em atoms} are the
elements covering $\hat{0}_P$). We denote by $\mathfrak{I}(A)$ and
$\mathfrak{F}(A)$ the order ideal and filter of $P$ generated by
an antichain $A$, respectively. If $Q$ is a subposet of $P$ then
$\bmin Q$ denotes the set of minimal elements of $Q$.

We call the empty antichain in $P$ and the one-element antichain
$\{\hat{0}_P\}$ the {\em trivial antichains\/} in $P$ because they
play in our study a role analogous to that played by the trivial
clutters in the theory of blocking sets.

We now recall some poset-theoretic blocker constructions. Let $j$
be a nonnegative integer less than $|P^{\mathrm a}|$. Given a
nontrivial antichain $A$ in $P$, define the antichain
\begin{equation}
\label{eq:1} \boldsymbol{\mathfrak{b}}_j(A):=\bmin\bigl\{b\in P:\
|\mathfrak{I}(b)\cap\mathfrak{I}(a)\cap P^{\mathrm a}|>j\ \
\forall a\in A\bigr\}\ .
\end{equation}
If $A$ is a trivial antichain in $P$ then the antichain
$\boldsymbol{\mathfrak{b}}_j(A)$ by definition is the other
trivial antichain.

The antichains $\boldsymbol{\mathfrak{b}}_j(A)$, defined
by~(\ref{eq:1}), serve as a poset-theoretic generalization of the
notion of set-theoretic blocker of a nontrivial clutter,
see~Matveev, 2003. From this point of view, the antichain
\begin{equation}\label{eq:2}
\boldsymbol{\mathfrak{b}}(A):=\boldsymbol{\mathfrak{b}}_0(A)
\end{equation}
bears a strong resemblance to its set-theoretic predecessor,
see~Bj\"{o}rner and Hultman, 2004, and Matveev, 2001.
Antichains~(\ref{eq:1}) admit a nice ordering, and some of the
structural and combinatorial properties of blockers~(\ref{eq:2})
in the Boolean lattices are clarified, see Remark~\ref{rem:1}.

The posets for which
\begin{equation*}
\boldsymbol{\mathfrak{b}}\bigl(\boldsymbol{\mathfrak{b}}(A)\bigr)=
A\ ,
\end{equation*}
for all antichains $A$, are characterized in~Bj\"{o}rner and
Hultman, 2004.

When we deal with construction~(\ref{eq:1}) related to a nontrivial antichain $A$, we are interested in the
nonemptiness and the cardinalities of the intersections $\mathfrak{I}(b)\cap\mathfrak{I}(a)\cap P^{\mathrm
a}$, for $b\in P-\{\hat{0}_P\}$ and $a\in A$, while the cardinalities of the sets $\mathfrak{I}(b)\cap
P^{\mathrm a}$ do not matter. To distinguish the objects we mainly study in the present paper from those
similar to~(\ref{eq:1}), we say that the antichain $\boldsymbol{\mathfrak{b}}_j(A)$ is an example of an {\em
absolute\/} poset-theoretic {\em $j$-blocker}; a more general definition is given in
Section~\ref{SectionOnConvexSubposetsAndAbsoluteBlockers}. Let $r$ be a rational number such that $0\leq
r<1$. A relative counterpart of $\boldsymbol{\mathfrak{b}}_j(A)$ is the antichain
\begin{equation}\label{eq:3}
\bmin\left\{b\in P-\{\hat{0}_P\}:\
\frac{|\mathfrak{I}(b)\cap\mathfrak{I}(a)\cap P^{\mathrm
a}|}{|\mathfrak{I}(b)\cap P^{\mathrm a}|}>r\ \ \forall a\in
A\right\}\ ;
\end{equation}
similar constructions form the subject of the present paper.

The study of poset-theoretic generalizations of set-theoretic
committees, undertaken in the paper, has been partly motivated by
the need for a more detailed analysis of building blocks of
decision rules in applied contradictory problems of pattern
recognition. See Duda et al., 2001, on the setting of the pattern
recognition problem and various methods to solve it.

Consider a finite nonempty collection $\boldsymbol{\mathcal{H}}:=\{\pmb{H}_1,\ldots,\pmb{H}_m\}$ of
codimension one linear subspaces $\pmb{H}_i:=\{\pmb{x}\in\mathbb{R}^n:\ \langle\pmb{p}_i,\pmb{x}\rangle=0\}$
in the {\em feature space\/} $\mathbb{R}^n$ with $n\geq 2$, where any two vectors from the rank $n$ set
$\{\pmb{p}_i:\ 1\leq i\leq m\}\subset\mathbb{R}^n$ are linearly independent;
$\langle\pmb{p}_i,\pmb{x}\rangle:=\sum_{j=1}^n p_{ij} x_j$. The connected components of the complement
$\mathbb{R}^n-\bigcup_{1\leq i\leq m}\pmb{H}_i$ of the {\em hyperplane arrangement\/}
$\boldsymbol{\mathcal{H}}$ are called the {\em regions\/} (or {\em chambers}) of $\boldsymbol{\mathcal{H}}$,
see e.g., Orlik and Terao, 1992.

We call the arrangement of oriented hyperplanes
$\boldsymbol{\mathcal{H}}$ (that is the set
$\boldsymbol{\mathcal{H}}$ for every hyperplane $\pmb{H}$ of which
``positive'' and ``negative sides'' of $\pmb{H}$ are
distinguished) a {\em training set}, if a partition
$\boldsymbol{\mathcal{H}}=\boldsymbol{\mathcal{A}}\dot\cup
\boldsymbol{\mathcal{B}}$ of $\boldsymbol{\mathcal{H}}$ into two
nonempty {\em training samples\/} $\boldsymbol{\mathcal{A}}$ and
$\boldsymbol{\mathcal{B}}$ is given. The hyperplanes from
$\boldsymbol{\mathcal{H}}$ are called the {\em training patterns}.
The training samples $\boldsymbol{\mathcal{A}}$ and
$\boldsymbol{\mathcal{B}}$ are thought of as subsets of two
disjoint {\em classes\/} $\mathbf{A}$ and $\mathbf{B}$,
respectively; these classes, in general, are sets of unknown
nature. We say that a pattern $\pmb{H}$ {\sl a priori\/} belongs
to the class $\mathbf{A}$ and it has the corresponding label
$\lambda(\pmb{H}):=-$, if $\pmb{H}\in\boldsymbol{\mathcal{A}}$;
the pattern $\pmb{H}$ {\sl a priori\/} belongs to the class
$\mathbf{B}$ and it has the label $\lambda(\pmb{H}):=+$, if
$\pmb{H}\in\boldsymbol{\mathcal{B}}$.

A region $\pmb{T}$ of $\boldsymbol{\mathcal{H}}$ lies on the {\em
positive side\/} of a hyperplane $\pmb{H}_i$, if the value
$\langle\pmb{e}_i,\pmb{v}\rangle$ is positive for some vector
$\pmb{v}\in\pmb{T}$, where the vector $\pmb{e}_i$ is defined by
$\pmb{e}_i:=-\pmb{p}_i$ for
$\pmb{H}_i\in\boldsymbol{\mathcal{A}}$, and by
$\pmb{e}_i:=\pmb{p}_i$ for $\pmb{H}_i\in\boldsymbol{\mathcal{B}}$.
Denote by $\boldsymbol{\mathcal{T}}_i^+$ the set of all regions
lying on the positive side of $\pmb{H}_i$.

We say that a subset of regions
$\boldsymbol{\mathcal{K}}^{\ast}:=\{\pmb{R}_1,\ldots,\pmb{R}_t\}$
of $\boldsymbol{\mathcal{H}}$ is a {\em committee\/} for
$\boldsymbol{\mathcal{H}}$ if for every $i$, $1\leq i\leq m$, it
holds
$|\boldsymbol{\mathcal{K}}^{\ast}\cap\boldsymbol{\mathcal{T}}_i^+|>
\tfrac{1}{2}|\boldsymbol{\mathcal{K}}^{\ast}|$. In this case a
system of representatives $\{\pmb{w}_k\in\pmb{R}_k:\ 1\leq k\leq
t\}$ is called a {\em committee\/} for the homogeneous system of
strict linear inequalities
$\{\langle\pmb{e}_i,\mathbf{x}\rangle>0:\ 1\leq i\leq m\}$.

Committees for such inequality systems were apparently first
introduced in Ablow and Kaylor, 1965, where it was proved that
such very useful collective generalizations of the notion of
solution do exist. Those notes laid the foundation of a branch of
the theory of pattern recognition; some of the surveys in the
committee mathematical methods and their applications are Khachai,
2004, Khachai et al., 2002, Mazurov, 1990, Mazurov et al., 1989,
and Mazurov and Khachai, 1999, 2004.

The {\em decision rule\/} $\mathfrak{r}$ is the mapping $\boldsymbol{\mathcal{H}}\to\{-,+\}$ under which
$\mathfrak{r}:\pmb{H}\mapsto\lambda(\pmb{H})$; in other words, such a rule must correctly recognize the
patterns from the training set.

Given a committee $\{\pmb{w}_k:\ 1\leq k\leq t\}$ for the
inequality system $\{\langle\pmb{e}_i,\mathbf{x}\rangle>0:\ 1\leq
i\leq m\}$, one defines the corresponding {\em committee decision
rule\/} $\mathfrak{r}$ in the following way: if $|\{\pmb{w}_k:\
\langle\pmb{p}_i,\pmb{w}_k\rangle>0\}|<\tfrac{t}{2}$ then
$\mathfrak{r}:\pmb{H}_i\mapsto -$; otherwise,
$\mathfrak{r}:\pmb{H}_i\mapsto +$.

When a new pattern, that is a new oriented hyperplane $\pmb{G}$,
is added to the training set $\boldsymbol{\mathcal{H}}$, the
domain and range of the decision rule $\mathfrak{r}$, associated
to the committee $\{\pmb{w}_k:\ 1\leq k\leq t\}$, extend over the
sets $\boldsymbol{\mathcal{H}}\dot\cup\pmb{G}$ and $\{-,0,+\}$,
respectively. The image of $\pmb{G}$ under $\mathfrak{r}$ is
determined depending on whether a majority of the vectors from
$\{\pmb{w}_k:\ 1\leq k\leq t\}$ lies on the positive side of
$\pmb{G}$. The case $\mathfrak{r}(\pmb{G})=0$ means that the new
pattern $\pmb{G}$ is not recognized.

In order to analyze the structural and combinatorial properties of
the family of all possible committees for the hyperplane
arrangement $\boldsymbol{\mathcal{H}}$ in detail, presumably, one
may consider the Boolean lattice $P$ of all subsets of the set of
regions of $\boldsymbol{\mathcal{H}}$. The language of the theory
of {\em oriented matroids\/} (which, for example, translates the
regions of $\boldsymbol{\mathcal{H}}$ to the {\em maximal
covectors\/} of a {\em realizable\/} oriented matroid) may be of
use; see Bj\"{o}rner et al., 1993, on oriented matroids. Recall
that the means of computing the rank of $P$, that is the number of
regions of $\boldsymbol{\mathcal{H}}$, are well-known (Zaslavsky,
1975). Nonempty subsets of regions, regarded as elements $b$ of
$P$, are committees for $\boldsymbol{\mathcal{H}}$ if and only if
the inequalities
\begin{equation*}
\frac{|\mathfrak{I}(b)\cap\mathfrak{I}(a)\cap P^{\mathrm
a}|}{|\mathfrak{I}(b)\cap P^{\mathrm a}|}>r
\end{equation*}
hold for all elements $a$ of the antichain
$A:=\{\boldsymbol{\mathcal{T}}_1^+,\ldots,\boldsymbol{\mathcal{T}}_m^+\}$
in $P$, under $r:=\tfrac{1}{2}$. From this point of view, the
elements of antichain~(\ref{eq:3}) are committees (which are
inclusion-minimal) of ``high quality'' for the arrangement
$\boldsymbol{\mathcal{H}}$.

In Section~\ref{SectionIntroducingRelativeBlockers} of this paper,
we introduce and discuss relative blocker constructions that
generalize constructions~(\ref{eq:3}). In
Section~\ref{SectionOnConvexSubposetsAndAbsoluteBlockers}, we turn
to their absolute predecessors going back to blocking sets and
set-theoretic blockers similar to~(\ref{eq:1}). In
Section~\ref{SectionOnBothConcepts}, we remark on a connection
between the concepts of absolute and relative blocking in posets.
In Section~\ref{sectiononproperties}, we analyze the structure of
relative blocker constructions, and we touch on the subject of
enumeration. Our exploration leads us to sequences of irreducible
fractions associated to the principal order ideals in posets which
are considered in Section~\ref{SectionOnFareySubsequences} and
studied, in the Boolean context, in
Section~\ref{SectionOnFareySubsequencesInBooleanContext}. It turns
out that all the familiar properties of the classical Farey
sequences of the theory of numbers are inherited by subsequences
of irreducible fractions whose nature is largely poset-theoretic.
In Section~\ref{SectionOnApplicationOfFareySubsequences}, we apply
Farey subsequences to relative blocker constructions in graded
posets.

If $Q$ is a subposet of $P$ then, throughout the paper, $\bmax Q$
stands for the set of maximal elements of $Q$. We denote by
$\mathfrak{A}_{\vartriangle}(P)$ and
$\mathfrak{A}_{\triangledown}(P)$ distributive lattices of all
antichains in $P$ defined in the following way. If $A'$ and $A''$
are antichains in $P$ then we set $A'\leq A''$ in
$\mathfrak{A}_{\vartriangle}(P)$ if and only if it holds
$\mathfrak{I}(A')\subseteq\mathfrak{I}(A'')$, and we set $A'\leq
A''$ in $\mathfrak{A}_{\triangledown}(P)$ if and only if it holds
$\mathfrak{F}(A')\subseteq\mathfrak{F}(A'')$. We use the notations
$\hat{0}_{\mathfrak{A}_{\vartriangle}(P)}$ and
$\hat{0}_{\mathfrak{A}_{\triangledown}(P)}$ to denote the least
elements of $\mathfrak{A}_{\vartriangle}(P)$ and
$\mathfrak{A}_{\triangledown}(P)$, respectively; we use the
similar notations $\hat{1}_{\mathfrak{A}_{\vartriangle}(P)}$ and
$\hat{1}_{\mathfrak{A}_{\triangledown}(P)}$ to denote the greatest
elements. The operations of meet in
$\mathfrak{A}_{\vartriangle}(P)$ and
$\mathfrak{A}_{\triangledown}(P)$ are denoted by
$\wedge_{\vartriangle}$ and $\wedge_{\triangledown}$,
respectively; in a similar manner, $\vee_{\vartriangle}$ and
$\vee_{\triangledown}$ stand for the operations of join. If $A'$
and $A''$ are antichains in $P$, then we have
$A'\wedge_{\vartriangle}A''=\bmax(\mathfrak{I}(A')\cap\mathfrak{I}(A''))$,
$A'\vee_{\vartriangle}A''=\bmax(A'\cup A'')$ and, in the dual
manner,
$A'\wedge_{\triangledown}A''=\bmin(\mathfrak{F}(A')\cap\mathfrak{F}(A''))$,
$A'\vee_{\triangledown}A''=\bmin(A'\cup A'')$.

Recall that in the present paper the least and greatest elements
of the lattice $\mathfrak{A}_{\triangledown}(P)$ are called the
trivial antichains in $P$;
$\hat{0}_{\mathfrak{A}_{\triangledown}(P)}$ is the empty antichain
in $P$, and $\hat{1}_{\mathfrak{A}_{\triangledown}(P)}$ is the
one-element antichain $\{\hat{0}_P\}$.

$\mathbb{Q}$ denotes rational numbers; $\mathbb{N}$, $\mathbb{P}$,
and $\mathbb{Z}$ stand for nonnegative, positive, and all
integers, respectively. $i|j$ means that an integer $i$ divides an
integer $j$; $i\bot j$ means that $i$ and $j$ are relatively
prime, and $\gcd(i,j)$ denotes the greatest common divisor of $i$
and $j$.

If $i$ and $j$ are positive integers then we denote by $\lbrack
i,j\rbrack$ the set $\{i,i+1,\ldots,j\}$.

If the poset $P$ is graded, with the rank function
$\rho:P\to\mathbb{N}$, then we write $\rho(P)$ instead of
$\rho(\hat{1}_P)$; further, given $j\in\{0\}\cup\lbrack
1,\rho(P)\rbrack$, we denote by $P^{(j)}$ the subset $\{p\in P:\
\rho(p)=j\}$. The layer $P^{(1)}=:P^{\mathrm a}$ is the set of
atoms of $P$.

Recall that a subposet $C$ of the poset $P$ is called {\em
convex\/} if the implication $x,z\in C$, $y\in P$, $x\leq y\leq z$
in $P$ $\Longrightarrow$ $y\in C$ holds for all elements $x,y,z\in
P$. We regard the empty subposet as a convex one.

The {\em M\"{o}bius function} (see, e.g.,~Chapter~IV of Aigner,
1979, Bj\"{o}rner et al., 1997, Greene, 1982, and Chapter~3 of
Stanley, 1997) $\mu_{P}:P\times P\to\mathbb{Z}$ is defined in the
following way: $\mu_P(x,x):=1$, for any $x\in P$; further, if
$z\in P$ and $x< z$ in $P$, then $\mu_P(x,z):=-\sum_{y\in P:\
x\leq y<z}\mu_P(x,y)$; finally, if $x\nleq z$ in $P$, then
$\mu_P(x,z):=0$.

We denote by $\mathbb{B}(n)$ the Boolean lattice of finite rank
$n\geq 1$. $\mathbb{V}_q(n)$ stands for the lattice of all
subspaces of a vector space of finite dimension $n\geq 1$ over a
finite field of $q$ elements. $\tbinom{j}{i}$ and
$\tbinom{j}{i}_{\! q}$ denote a binomial and $q$-binomial
coefficient, respectively.

Finally, $r$ always denotes a rational number such that $0\leq
r<1$.

\section{Relative $r$-blockers}
\label{SectionIntroducingRelativeBlockers}

Let
\begin{equation}
\label{eq:4}
\omega:\mathfrak{A}_{\vartriangle}(P)\to\{-1\}\cup\mathbb{N}
\end{equation}
be a map such that
\begin{equation}
\label{eq:5} \hat{0}_{\mathfrak{A}_{\vartriangle}(P)}\mapsto -1\
,\ \ \ \{\hat{0}_P\}\mapsto 0\ ;
\end{equation}
and for any antichains $A'$ and $A''$ in $P$ such that
$\{\hat{0}_P\}<A'\leq A''$ in $\mathfrak{A}_{\vartriangle}(P)$, it
holds
\begin{equation}
\label{eq:6} 0<\omega(A') \leq \omega(A'')\ .
\end{equation}

From now on, $\omega$ always means map~(\ref{eq:4}) satisfying
constraints~(\ref{eq:5}) and (\ref{eq:6}). Some relevant examples
of $\omega$ follow:
\begin{itemize}
\item
\begin{equation*}
\omega: A\mapsto\rho_{\vartriangle}(A)-1=|\mathfrak{I}(A)|-1\ ,
\end{equation*}
where $\rho_{\vartriangle}(A)$ denotes the rank of an element $A$
in the lattice $\mathfrak{A}_{\vartriangle}(P)$;
\item
\begin{equation}
\label{eq:7} \omega: A\mapsto
\rho_{\vartriangle}(A\wedge_{\vartriangle}P^{\mathrm
a})-1=\begin{cases}-1,&\text{if
$A=\hat{0}_{\mathfrak{A}_{\vartriangle}(P)}$}\ ,\\
|\mathfrak{I}(A)\cap P^{\mathrm a}|,&\text{if
$A\neq\hat{0}_{\mathfrak{A}_{\vartriangle}(P)}$}\ ;
\end{cases}
\end{equation}
\item
\begin{equation}\label{eq:8} \omega:A\mapsto\begin{cases}-1,& \text{if
$A=\hat{0}_{\mathfrak{A}_{\vartriangle}(P)}$}\ ,\\ \max_{a\in
A}\rho(a),& \text{if
$A\neq\hat{0}_{\mathfrak{A}_{\vartriangle}(P)}$}\ ,\end{cases}
\end{equation}
if $P$ is graded, with the rank function $\rho$.
\end{itemize}

The maps $\omega$ defined by~{\rm(\ref{eq:4})-(\ref{eq:6})} are
sometimes well expressed in terms of {\em incidence functions};
see, e.g.,~Chapter~IV of Aigner, 1979, and Chapter~3 of Stanley,
1997, on incidence functions of posets.

Throughout the paper, we write $\rho$ instead of $\omega$ when we
deal exclusively with map~(\ref{eq:4}) defined by~(\ref{eq:8}). If
$\{a\}$ is a one-element antichain in $P$ then we write
$\omega(a)$ instead of $\omega(\{a\})$, and we write $\omega(P)$
instead of
$\omega(\hat{1}_P)=\omega(\hat{1}_{\mathfrak{A}_{\vartriangle}(P)})$.

\begin{defn}
\label{defn:1} Let $A$ be a subset of $P$.
\begin{itemize}
\item[\rm (i)]
If $A$ is nonempty and $A\neq\{\hat{0}_P\}$, then an element $b\in
P-\{\hat{0}_P\}$ is a {\em relatively $r$-blocking element for $A$
in $P$} {\rm(}w.r.t. a map $\omega${\rm)} if, for every $a\in
A-\{\hat{0}_P\}$, it holds
\begin{equation}
\label{eq:9}
\frac{\omega(\{b\}\wedge_{\vartriangle}\{a\})}{\omega(b)}>r\ .
\end{equation}
\item[\rm (ii)]
If $A=\{\hat{0}_P\}$ then $A$ has no {\em relatively $r$-blocking
elements in $P$}.
\item[\rm (iii)]
If $A$ is empty then every element of $P$ is a {\em relatively
$r$-blocking element for $A$ in $P$}.
\end{itemize}
\end{defn}

{\small
\begin{rem} \label{rem:4}Let $A$ be a nonempty subset of\/
$\mathbb{B}(n)-\{\hat{0}_{\mathbb{B}(n)}\}$. An element\/ $b\in
\mathbb{B}(n)-\{\hat{0}_{\mathbb{B}(n)}\}$ is a relatively
$r$-blocking element for $A$ in $\mathbb{B}(n)$, w.r.t. either of
the maps $\omega$ defined by~{\rm(\ref{eq:7})}
and~{\rm(\ref{eq:8})}, if and only if the set
$\mathfrak{I}(b)\cap\mathbb{B}(n)^{(1)}$ is an $r$-committee for
the family $\{\mathfrak{I}(a)\cap\mathbb{B}(n)^{(1)}:\ a\in A\}$,
that is, it holds
\begin{equation*}
|\mathfrak{I}(b)\cap\mathfrak{I}(a)\cap\mathbb{B}(n)^{(1)}|>r
\cdot|\mathfrak{I}(b)\cap\mathbb{B}(n)^{(1)}|\ ,
\end{equation*}
for all $a\in A$.
\end{rem}
}

We denote the subposet of $P$ consisting of all relatively
$r$-blocking elements for $A$, w.r.t. a map $\omega$, by
$\mathbf{I}_r(P,A;\omega)$. Given $a\in P$, we write
$\mathbf{I}_r(P,a;\omega)$ instead of
$\mathbf{I}_r(P,\{a\};\omega)$. If $k\in\lbrack
1,\omega(P)\rbrack$ then we denote by
$\mathbf{I}_{r,k}(P,A;\omega)$ the subposet
$\{b\in\mathbf{I}_r(P,A;\omega):\ \omega(b)=k\}$.

If $A$ is a nonempty subset of $P-\{\hat{0}_P\}$ then
Definition~\ref{defn:1} implies
$\mathbf{I}_r(P,A;\omega)=\mathbf{I}_r(P,\bmin A;\omega)$; this is
the reason why we are primarily interested in relatively
$r$-blocking elements for antichains.

If $A$ is a nontrivial antichain in $\mathbb{B}(n)$ then its order
ideal $\mathfrak{I}(A)$ is assigned the isomorphic {\em face
poset\/} of the {\em abstract simplicial complex} whose {\em
facets\/} are the sets from the family
$\{\mathfrak{I}(a)\cap\mathbb{B}(n)^{(1)}:\ a\in A\}$. See,
e.g.,~Billera and Bj\"{o}rner, 1997, Bj\"{o}rner, 1995, Bruns and
Herzog, 1998, Buchstaber and Panov, 2004, Hibi, 1992, Miller and
Sturmfels, 2004, Stanley, 1996, and Ziegler, 1998, on simplicial
complexes.

The following proposition lists some observations.

\begin{prop}
\label{prop:1}
\begin{itemize}
\item[\rm (i)] If $A$ is a nontrivial antichain in $P$, then it holds
\begin{equation*}
\mathbf{I}_r(P,A;\omega)=\bigcap_{a\in A}\mathbf{I}_r(P,a;\omega)\
,
\end{equation*}
for any map $\omega$.
\item[\rm (ii)]
If $A'$ and $A''$ are antichains in $P$ and $A'\leq A''$ in
$\mathfrak{A}_{\triangledown}(P)$, then
$\mathbf{I}_r(P,A';\omega)\supseteq\mathbf{I}_r(P,A'';\omega)$,
for any map $\omega$.
\item[\rm (iii)]
Let $r',r''\in\mathbb{Q}$, $0\leq r'\leq r''<1$. For any antichain
$A$ in $P$, and for any map $\omega$, it holds
$\mathbf{I}_{r'}(P,A;\omega)\supseteq\mathbf{I}_{r''}(P,A;\omega)$.
\end{itemize}
\end{prop}

The minimal elements of the subposets $\mathbf{I}_r(P,A;\omega)$
of the poset $P$ are of interest.

\begin{defn}
\label{defn:2}
\begin{itemize}
\item[\rm(i)]
The {\em relative $r$-blocker map on
$\mathfrak{A}_{\triangledown}(P)$} {\rm(}w.r.t. a map
$\omega${\rm)} is the map $\boldsymbol{\mathfrak{y}}_{\mathit
r}:\mathfrak{A}_{\triangledown}(P)\to\mathfrak{A}_{\triangledown}(P)$,
defined by
\begin{multline*}
A\mapsto\bmin \mathbf{I}_r(P,A;\omega)\\=\bmin \left\{b\in
P-\{\hat{0}_P\}:\
\frac{\omega(\{b\}\wedge_{\vartriangle}\{a\})}{\omega(b)}>r \ \
\forall a\in A\right\}
\end{multline*}
if $A$ is nontrivial, and
\begin{equation*}
\hat{0}_{\mathfrak{A}_{\triangledown}(P)}\mapsto
\hat{1}_{\mathfrak{A}_{\triangledown}(P)}\ ,\ \ \
\hat{1}_{\mathfrak{A}_{\triangledown}(P)}\mapsto
\hat{0}_{\mathfrak{A}_{\triangledown}(P)}\ .
\end{equation*}
\item[\rm(ii)]
Given an antichain $A$ in $P$, the antichain
$\boldsymbol{\mathfrak{y}}_{\mathit r}(A)$ is called the {\em
relative $r$-blocker\/} {\rm(}w.r.t. the map $\omega${\rm)} {\em
of $A$ in $P$}; the elements of
$\boldsymbol{\mathfrak{y}}_{\mathit r}(A)$ are called the {\em
minimal relatively $r$-blocking elements\/} {\rm(}w.r.t. the map
$\omega${\rm)} {\em for $A$ in $P$}.
\end{itemize}
\end{defn}

In addition to the minimal relatively $r$-blocking elements, the
relatively $r$-blocking elements $b$ for $A$ in $P$ with the
minimum value of $\omega(b)$ can be of particular interest.

The following statement is a consequence of
Proposition~\ref{prop:1}(ii,iii). It particularly states that the
relative $r$-blocker map is order-reversing.
\begin{cor}
Let $r',r''\in\mathbb{Q}$, $0\leq r'\leq r''<1$. Let $A'$ and
$A''$ be antichains in $P$ such that $A'\leq A''$ in
$\mathfrak{A}_{\triangledown}(P)$.  The relation
\begin{equation*}
\boldsymbol{\mathfrak{y}}_{r''}(A'')\leq
\boldsymbol{\mathfrak{y}}_{r'}(A'')\leq\boldsymbol{\mathfrak{y}}_{r'}(A')
\end{equation*}
holds in $\mathfrak{A}_{\triangledown}(P)$.
\end{cor}

Let $A$ be a nontrivial antichain in $P$. If the relative
$r$-blocker $\boldsymbol{\mathfrak{y}}_{\mathit r}(A)$ of $A$ in
$P$ (w.r.t. a map $\omega$) is not
$\hat{0}_{\mathfrak{A}_{\triangledown}(P)}$, then $A$ is a subset
of relatively $r'$-blocking elements for the antichain
$\boldsymbol{\mathfrak{y}}_{\mathit r}(A)$, for some
$r'\in\mathbb{Q}$. Indeed, for each $a\in A$ and for all
$b\in\boldsymbol{\mathfrak{y}}_{\mathit r}(A)$, we by~(\ref{eq:9})
have
\begin{equation*}
\frac{\omega(\{a\}\wedge_{\vartriangle}\{b\})}{\omega(a)}>
r\cdot\frac{\omega(b)}{\omega(a)}\geq
r\cdot\frac{\min_{p\in\boldsymbol{\mathfrak{y}}_{\mathit
r}(A)}\omega(p)}{\max_{p\in A}\omega(p)}\ ,
\end{equation*}
and this observation implies the following statement.

\begin{prop} If $A$ is a nontrivial antichain in $P$ and
$\boldsymbol{\mathfrak{y}}_{\mathit
r}(A)\neq\hat{0}_{\mathfrak{A}_{\triangledown}(P)}$, w.r.t. a map
$\omega$, then
\begin{equation*}
A\subseteq\mathbf{I}_{r'}\bigl(P,\boldsymbol{\mathfrak{y}}_{\mathit
r}(A);\omega\bigr)\ ,
\end{equation*}
where
$r':=r\cdot\tfrac{\min_{p\in\boldsymbol{\mathfrak{y}}_{\mathit
r}(A)}\omega(p)}{\max_{p\in A}\omega(p)}$.
\end{prop}

\section{Absolute $j$-blockers and convex subposets}
\label{SectionOnConvexSubposetsAndAbsoluteBlockers}

Let $A$ be a nontrivial antichain in $P$. Let $h$ and $k$ be
positive integers such that $h\leq k\leq\omega(P)$, for some map
$\omega$. In the following sections of the paper we will make use
of the auxiliary subposet
\begin{equation}
\label{eq:10} \left\{b\in P:\ \omega(b)=k,\
\omega\bigl(\{b\}\wedge_{\vartriangle}\{a\}\bigr)\geq h\ \ \forall
a\in A\right\}\ .
\end{equation}
We can consider this subposet, in an equivalent way, as the
intersection
\begin{multline}
\label{eq:11} \Bigl(\ \bigl\{b\in P:\ \omega(b)>k-1\bigr\}\ -\
\bigl\{b\in P:\ \omega(b)>k\bigr\}\ \Bigr)\\  \cap \bigl\{b\in P:\
\omega\bigl(\{b\}\wedge_{\vartriangle}\{a\}\bigr)> h-1\ \ \forall
a\in A\bigr\}\ .
\end{multline}
Each component of expression~(\ref{eq:11}) can be described in
terms of {\em absolute blocking}. Indeed, given a nontrivial
antichain $A$ in $P$ and a nonnegative integer $j$ less than
$\omega(P)$, define the {\em absolute $j$-blocker\/} (w.r.t. the
map $\omega$) {\em of $A$ in $P$}, denoted by
$\boldsymbol{\mathsf{b}}_j(A)$, in the following way:
\begin{equation}
\label{eq:12} \boldsymbol{\mathsf{b}}_j(A):=\bmin\bigl\{b\in P:\
\omega\bigl(\{b\}\wedge_{\vartriangle}\{a\}\bigr)>j\ \ \forall
a\in A\bigr\}\ .
\end{equation}
For any element $b\in\mathfrak{F}\bigl(\boldsymbol{\mathsf{b}}_j(A)\bigr)$, we have
$\omega\bigl(\{b\}\wedge_{\vartriangle}\{a\}\bigr)>j$, for all $a\in A$. A particular example of absolute
$j$-blocker~(\ref{eq:12}) is the construction defined by~(\ref{eq:1}) and implicitly involving the map
$\omega$ defined by~(\ref{eq:7}). We set $\boldsymbol{\mathsf{b}}_{\omega(P)}(A):=
\hat{0}_{\mathfrak{A}_{\triangledown}(P)}$. Note that
\begin{equation}
\boldsymbol{\mathsf{b}}_j(A)=\bigwedge_{a\in
A}\boldsymbol{\mathsf{b}}_j(a)
\end{equation}
in $\mathfrak{A}_{\triangledown}(P)$; we write
$\boldsymbol{\mathsf{b}}_j(a)$ instead of
$\boldsymbol{\mathsf{b}}_j(\{a\})$.

If the trivial antichains in $P$ must be taken into consideration
then we set
\begin{equation}
\label{eq:14}
\boldsymbol{\mathsf{b}}_j\!\left(\hat{0}_{\mathfrak{A}_{\triangledown}(P)}\right)
:=\hat{1}_{\mathfrak{A}_{\triangledown}(P)}\ ,\ \ \ \
\boldsymbol{\mathsf{b}}_j\!\left(\hat{1}_{\mathfrak{A}_{\triangledown}(P)}\right)
:=\hat{0}_{\mathfrak{A}_{\triangledown}(P)}\ .
\end{equation}

Given an antichain $A$ in $P$ and a map $\omega$, we call the
elements of the order filter
$\mathfrak{F}\bigl(\boldsymbol{\mathsf{b}}_j(A)\bigr)$ the {\em
absolutely $j$-blocking elements for $A$ in $P$} (w.r.t. the map
$\omega$). The elements of the order filter
$\mathfrak{F}\bigl(\boldsymbol{\mathfrak{b}}(A)\bigr)$, where the
antichain $\boldsymbol{\mathfrak{b}}(A)$ is defined
by~(\ref{eq:2}), were called in~Matveev, 2001, the {\em
intersecters for $A$ in $P$}.

If $P$ is graded, and if the map $\omega$ is defined
by~(\ref{eq:8}) then, given a nontrivial one-element antichain
$\{a\}$ in $P$, we have
\begin{equation*}
\boldsymbol{\mathsf{b}}_j(a)=\mathfrak{I}(a)\cap P^{(j+1)}\ .
\end{equation*}

The {\em absolute $j$-blocker map\/} $\boldsymbol{\mathsf{b}}_j:\
\mathfrak{A}_{\triangledown}(P)\to\mathfrak{A}_{\triangledown}(P)$
is order-reversing, w.r.t. any map $\omega$. If $A$ is an
arbitrary antichain in $P$ then for any nonnegative integers $i$
and $j$ such that $i\leq j<\omega(P)$, the relation
\begin{equation}\label{eq:15}
\boldsymbol{\mathsf{b}}_i(A)\geq \boldsymbol{\mathsf{b}}_j(A)
\end{equation}
holds in $\mathfrak{A}_{\triangledown}(P)$.

If $A$ is a trivial antichain in $P$ then convention~(\ref{eq:14})
implies
$\boldsymbol{\mathsf{b}}_j\bigl(\boldsymbol{\mathsf{b}}_j(A)\bigr)$
$=A$. Now, let $A$ be a nontrivial antichain. If
$\boldsymbol{\mathsf{b}}_j(A)=\hat{0}_{\mathfrak{A}_{\triangledown}(P)}$,
then we have
$\boldsymbol{\mathsf{b}}_j\bigl(\boldsymbol{\mathsf{b}}_j(A)\bigr)=
\hat{1}_{\mathfrak{A}_{\triangledown}(P)}>A$ in
$\mathfrak{A}_{\triangledown}(P)$. Finally, suppose that
$\boldsymbol{\mathsf{b}}_j(A)$ is a nontrivial antichain in $P$.
On the one hand, for each $a\in A$ and for all
$b\in\boldsymbol{\mathsf{b}}_j(A)$, we have
$\omega\bigl(\{a\}\wedge_{\vartriangle}\{b\}\bigr)>j$. On the
other hand, (\ref{eq:12}) implies
\begin{equation}\label{eq:16}
\boldsymbol{\mathsf{b}}_j\bigl(\boldsymbol{\mathsf{b}}_j(A)\bigr)=
\bmin\Bigl\{g\in P:\
\omega\bigl(\{g\}\wedge_{\vartriangle}\{b\}\bigr)>j\ \ \forall
b\in\boldsymbol{\mathsf{b}}_j(A)\Bigr\}\ .
\end{equation}
Hence we have
\begin{equation}
\label{eq:17}
\boldsymbol{\mathsf{b}}_j\bigl(\boldsymbol{\mathsf{b}}_j(A)\bigr)\geq
A
\end{equation}
in $\mathfrak{A}_{\triangledown}(P)$, for any
$A\in\mathfrak{A}_{\triangledown}(P)$.

Since $\boldsymbol{\mathsf{b}}_j$ is order-reversing and
(\ref{eq:17}) holds, the technique of the Galois correspondence
(see, e.g.,~Sections~IV.3.B,A of Aigner, 1979) can be applied to
the absolute $j$-blocker map $\boldsymbol{\mathsf{b}}_j$ on
$\mathfrak{A}_{\triangledown}(P)$:

\begin{prop} \label{prop:2}
Let $\boldsymbol{\mathsf{b}}_j:\mathfrak{A}_{\triangledown}(P)
\to\mathfrak{A}_{\triangledown}(P)$ be the absolute $j$-blocker
map on $\mathfrak{A}_{\triangledown}(P)$, w.r.t. a map $\omega$.
\begin{itemize}
\item[\rm(i)]
The composite map
$\boldsymbol{\mathsf{b}}_j\circ\boldsymbol{\mathsf{b}}_j$ is a
closure operator on $\mathfrak{A}_{\triangledown}(P)$.
\item[\rm(ii)]
The image
$\boldsymbol{\mathsf{b}}_j\left(\mathfrak{A}_{\triangledown}(P)\right)$
of the lattice $\mathfrak{A}_{\triangledown}(P)$ under the map
$\boldsymbol{\mathsf{b}}_j$ is a self-dual lattice; the
restriction of the map $\boldsymbol{\mathsf{b}}_j$ to
${\boldsymbol{\mathsf{b}}_j(\mathfrak{A}_{\triangledown}(P))}$ is
an anti-automorphism of
$\boldsymbol{\mathsf{b}}_j\bigl(\mathfrak{A}_{\triangledown}(P)\bigr)$.
As a consequence, for any antichain
$B\in\boldsymbol{\mathsf{b}}_j\left(\mathfrak{A}_{\triangledown}(P)\right)$
it holds
$\boldsymbol{\mathsf{b}}_j\bigl(\boldsymbol{\mathsf{b}}_j(B)\bigr)=B$.

The lattice
$\boldsymbol{\mathsf{b}}_j\bigl(\mathfrak{A}_{\triangledown}(P)\bigr)$
is a sub-meet-semilattice of $\mathfrak{A}_{\triangledown}(P)$.
\item[\rm(iii)]
For any
$B\in\boldsymbol{\mathsf{b}}_j\bigl(\mathfrak{A}_{\triangledown}(P)\bigr)$,
its preimage $(\boldsymbol{\mathsf{b}}_j)^{-1}(B)$ in
$\mathfrak{A}_{\triangledown}(P)$ under the map
$\boldsymbol{\mathsf{b}}_j$ is a convex sub-join-semilattice of
$\mathfrak{A}_{\triangledown}(P)$; the greatest element of
$(\boldsymbol{\mathsf{b}}_j)^{-1}(B)$ is
$\boldsymbol{\mathsf{b}}_j(B)$.
\end{itemize}
\end{prop}

\begin{proof} Assertions~(i) and~(ii) are consequences
of~Propositions~4.36 and~4.26 of Aigner, 1979.

To prove assertion~(iii), pick arbitrary elements
$A',A''\in(\boldsymbol{\mathsf{b}}_j)^{-1}(B)$, where
$B=\boldsymbol{\mathsf{b}}_j(A)$, for some
$A\in\mathfrak{A}_{\triangledown}(P)$, and note that
$\boldsymbol{\mathsf{b}}_j(A'\vee_{\triangledown}
A'')=\boldsymbol{\mathsf{b}}_j(A')
\wedge_{\triangledown}\boldsymbol{\mathsf{b}}_j(A'')=B$. Thus,
$(\boldsymbol{\mathsf{b}}_j)^{-1}(B)$ is a sub-join-semilattice of
$\mathfrak{A}_{\triangledown}(P)$. If
$B=\hat{0}_{\mathfrak{A}_{\triangledown}(P)}$ then
$\boldsymbol{\mathsf{b}}_j(B)=\hat{1}_{\mathfrak{A}_{\triangledown}(P)}$
is the greatest element of $(\boldsymbol{\mathsf{b}}_j)^{-1}(B)$.
If $B=\hat{1}_{\mathfrak{A}_{\triangledown}(P)}$ then
$(\boldsymbol{\mathsf{b}}_j)^{-1}(B)$ is the one-element subposet
$\{\hat{0}_{\mathfrak{A}_{\triangledown}(P)}\}\subset
\mathfrak{A}_{\triangledown}(P)$. Finally, if $B$ is a nontrivial
antichain in $P$ then the element
$\boldsymbol{\mathsf{b}}_j(B)=\boldsymbol{\mathsf{b}}_j\bigl(\boldsymbol{\mathsf{b}}_j(A)\bigr)$
is by~(\ref{eq:16}) the greatest element of
$(\boldsymbol{\mathsf{b}}_j)^{-1}(B)$. Since the map
$\boldsymbol{\mathsf{b}}_j$ is order-reversing, the subposet
$(\boldsymbol{\mathsf{b}}_j)^{-1}(B)$ of
$\mathfrak{A}_{\triangledown}(P)$ is convex.
\end{proof}

{\small
\begin{rem}\label{rem:1}
Let $A$ be an arbitrary antichain in the Boolean lattice
$\mathbb{B}(n)$. The antichain $\boldsymbol{\mathfrak{b}}(A)$
defined by~{\rm(\ref{eq:2})} satisfies the equality
$|\mathfrak{F}(A)|+|\mathfrak{F}(\boldsymbol{\mathfrak{b}}(A))|$
$=2^n$. As a consequence, we have $A=\boldsymbol{\mathfrak{b}}(A)$
if and only if it holds $|\mathfrak{F}(A)|=2^{n-1}$. In other
words, the layer
$\mathfrak{A}_{\triangledown}(\mathbb{B}(n))^{(2^{n-1})}$ of
$\mathfrak{A}_{\triangledown}(\mathbb{B}(n))$ is the set of {\em
fixed points} of the map $\boldsymbol{\mathfrak{b}}$. Indeed, we
have
$\boldsymbol{\mathfrak{b}}\bigl(\mathfrak{A}_{\triangledown}(\mathbb{B}(n))\bigr)=
\mathfrak{A}_{\triangledown}(\mathbb{B}(n))$, and our observations
follow immediately from Proposition~{\rm\ref{prop:2}(ii)}.
\end{rem}
}

We now return to consider poset~(\ref{eq:10}),(\ref{eq:11}). Note
that
\begin{equation*}
\begin{split}
\bigl\{b\in P:\
\omega(b)>k-1\bigr\}&=\mathfrak{F}(\boldsymbol{\mathsf{b}}_{k-1}(\hat{1}_P))\
,\\ \bigl\{b\in P:\
\omega(b)>k\bigr\}&=\mathfrak{F}(\boldsymbol{\mathsf{b}}_k(\hat{1}_P))\
,
\\
\bigl\{b\in P:\ \omega\bigl(\{b\}\wedge_{\vartriangle}\{a\}\bigr)>
h-1\ \ \forall a\in
A\bigr\}&=\mathfrak{F}\bigl(\boldsymbol{\mathsf{b}}_{h-1}(A)\bigr)\
;
\end{split}
\end{equation*}
therefore we obtain
\begin{multline}
\label{eq:18} \{b\in P:\ \omega(b)=k,\
\omega\bigl(\{b\}\wedge_{\vartriangle}\{a\}\bigr)\geq h\ \ \forall
a\in A\}\\ = \Bigl(
\mathfrak{F}\bigl(\boldsymbol{\mathsf{b}}_{k-1}(\hat{1}_P)\bigr)-
\mathfrak{F}\bigl(\boldsymbol{\mathsf{b}}_k(\hat{1}_P)\bigr)
\Bigr)\cap\mathfrak{F}\bigl(\boldsymbol{\mathsf{b}}_{h-1}(A)\bigr)\\
=
\mathfrak{F}\bigl(\boldsymbol{\mathsf{b}}_{k-1}(\hat{1}_P)\wedge_{\triangledown}\boldsymbol{\mathsf{b}}_{h-1}(A)\bigr)-
\mathfrak{F}\bigl(\boldsymbol{\mathsf{b}}_{k}(\hat{1}_P)\bigr)\ .
\end{multline}

Since $\boldsymbol{\mathsf{b}}_{k-1}(\hat{1}_P)\geq
\boldsymbol{\mathsf{b}}_{k}(\hat{1}_P)$ in
$\mathfrak{A}_{\triangledown}(P)$, by~(\ref{eq:15}), the second
line in expression~(\ref{eq:18}) describes an intersection of
convex subposets of $P$; hence the subposet presented in the first
line of~(\ref{eq:18}) is convex.

Again, let $h$ and $k$ be positive integers such that $h\leq
k\leq\omega(P)$. Let $\{a\}$ be a nontrivial one-element antichain
in $P$. In the following, in addition to
subposet~(\ref{eq:10}),(\ref{eq:11}),(\ref{eq:18}), we will also
need the convex subposet
\begin{multline}\label{eq:19}
\left\{b\in P:\ \omega(b)=k,\
\omega\bigl(\{b\}\wedge_{\vartriangle}\{a\}\bigr)=h\right\}\\=
\Bigl(\mathfrak{F}\bigl(\boldsymbol{\mathsf{b}}_{k-1}(\hat{1}_P)\bigr)-
\mathfrak{F}\bigl(\boldsymbol{\mathsf{b}}_{k}(\hat{1}_P)\bigr)\Bigr)\cap
\Bigl(\mathfrak{F}\bigl(\boldsymbol{\mathsf{b}}_{h-1}(a)\bigr)-
\mathfrak{F}\bigl(\boldsymbol{\mathsf{b}}_{h}(a)\bigr)\Bigr)
\\=
\mathfrak{F}\bigl(\boldsymbol{\mathsf{b}}_{k-1}(\hat{1}_P)
\wedge_{\triangledown}\boldsymbol{\mathsf{b}}_{h-1}(a)\bigr)-
\mathfrak{F}\bigl(\boldsymbol{\mathsf{b}}_{k}(\hat{1}_P)\vee_{\triangledown}
\boldsymbol{\mathsf{b}}_{h}(a)\bigr)\ .
\end{multline}

{\small
\begin{rem}\label{rem:2}
Let $h,k,m$ and $n$ be positive integers such that $m\leq n$ and
$h\leq k\leq n$. Recall that if $\{a\}$ is a nontrivial
one-element  antichain in $\mathbb{V}_q(n)$ with $\rho(a)=:m$,
then we have
\begin{equation*}
\bigl\{b\in\mathbb{V}_q(n):\ \rho(b)=k,\ \rho(b\wedge a)\geq
h\bigr\}=\mathfrak{F}\left(\mathfrak{I}(a)
\cap\mathbb{V}_q(n)^{(h)}\right)\cap\mathbb{V}_q(n)^{(k)}
\end{equation*}
and
\begin{equation*}
\left|\bigl\{b\in\mathbb{V}_q(n):\ \rho(b)=k,\ \rho(b\wedge a)\geq
h\bigr\}\right|=\sum_{j\in\lbrack h,k \rbrack}\binom{m}{j}_{\!\!
q}\binom{n-m}{k-j}_{\!\! q} q^{(m-j)(k-j)}\ .
\end{equation*}
Similarly, we have
\begin{multline*}
\bigl\{b\in\mathbb{V}_q(n):\ \rho(b)=k,\ \rho(b\wedge a)=h\bigr\}\\ = \Bigl(\mathfrak{F}\left(\mathfrak{I}(a)
\cap\mathbb{V}_q(n)^{(h)}\right)- \mathfrak{F}\left(\mathfrak{I}(a)\cap\mathbb{V}_q(n)^{(h+1)}\right)\Bigr)
\cap\mathbb{V}_q(n)^{(k)}
\end{multline*}
and
\begin{equation*}
\left|\bigl\{b\in\mathbb{V}_q(n):\ \rho(b)=k,\ \rho(b\wedge
a)=h\bigr\}\right|=\binom{m}{h}_{\!\! q}\binom{n-m}{k-h}_{\!\! q}
q^{(m-h)(k-h)}\ .
\end{equation*}
These expressions for the cardinalities of subposets have a direct
connection with the {\rm(}$q$-{\rm)}{\em Vandermonde's
convolution}, see, e.g.,~{\rm Section~4 of Andrews, 1974}.
\end{rem}
}

\section{Connection between concepts of absolute and relative blocking}
\label{SectionOnBothConcepts}

It follows from Definition~\ref{defn:2} that the relative
$0$-blocker $\boldsymbol{\mathfrak{y}}_0(A)$ of a nontrivial
antichain $A$ in $P$, w.r.t. an arbitrary map $\omega$, is nothing
else than the absolute $0$-blocker $\boldsymbol{\mathfrak{b}}(A)$
of $A$ in $P$, defined by~(\ref{eq:2}) and considered
in~Bj\"{o}rner et al., 2004, 2005, and Matveev, 2001. Moreover, if
$\boldsymbol{\mathfrak{b}}(A)\subseteq P^{\mathrm a}$ then
$\bigcap_{a\in
A}\mathfrak{I}(a)-\{\hat{0}_P\}\subseteq\mathbf{I}_r(P,A;\omega)$
and $\boldsymbol{\mathfrak{y}}_r(A)=\boldsymbol{\mathfrak{b}}(A)$,
for any value of the parameter $r$.

Again, let $A$ be a nontrivial antichain in $P$, and let
$j\in\mathbb{N}$, $j<\omega(P)$, for some map $\omega$. If
$\boldsymbol{\mathsf{b}}_j(A)\neq\hat{0}_{\mathfrak{A}_{\vartriangle}(P)}$
then, for all $b\in\boldsymbol{\mathsf{b}}_j(A)$ and for all $a\in
A$, we by~(\ref{eq:12}) have
\begin{equation*}
\begin{split}
\frac{\omega(\{b\}\wedge_{\vartriangle}\{a\})}{\omega(b)}&>
\frac{j}{\max_{p\in\boldsymbol{\mathsf{b}}_j(A)}\omega(p)}\ ,\\
\frac{\omega(\{a\}\wedge_{\vartriangle}\{b\})}{\omega(a)}&>
\frac{j}{\max_{p\in A}\omega(p)}\ ;
\end{split}
\end{equation*}
if $\boldsymbol{\mathfrak{y}}_{\mathit
r}(A)\neq\hat{0}_{\mathfrak{A}_{\vartriangle}(P)}$ then, for each
$b\in\boldsymbol{\mathfrak{y}}_{\mathit r}(A)$ and for all $a\in
A$, we by~(\ref{eq:9}) have
\begin{equation*}
\omega\bigl(\{b\}\wedge_{\vartriangle}\{a\}\bigr)>r\cdot\omega(b)\geq
r\cdot\min_{p\in\boldsymbol{\mathfrak{y}}_{\mathit
r}(A)}\omega(p)\ .
\end{equation*}

\section{Structure and enumeration}
\label{sectiononproperties}

We now turn to explore the structure of the subposets of
relatively $r$-blocking elements.

For $k\in\mathbb{P}$ such that $k\leq \omega(P)$, define the
integer
\begin{equation}\label{eq:20}
\nu(r\cdot k):=
\begin{cases}
\lceil r\cdot k\rceil&\text{\ \ if $r\cdot k\not\in\mathbb{N}$}\
,\\ r\cdot k+1,& \text{\ \ if $r\cdot k\in\mathbb{N}$}\ .
\end{cases}
\end{equation}
If $A$ is a nontrivial antichain in $P$, then it follows from
Definition~\ref{defn:1}(i) that it holds
\begin{multline}\label{eq:21}
\mathbf{I}_r(P,A;\omega)=\bigcup_{1\leq k\leq \omega(P)}\\
\bigcup_{\nu(r\cdot k)\leq h\leq\max_{a\in A}\omega(a)}
\bigl\{b\in P:\ \omega(b)=k,\
\omega\bigl(\{b\}\wedge_{\vartriangle}\{a\}\bigr)\geq h\ \ \forall
a\in A\bigr\}\ .
\end{multline}
Recall that for any values of $h$ and $k$ appearing in the above
expression, the structure of the poset $\bigl\{b\in P:\
\omega(b)=k,\
\omega\bigl(\{b\}\wedge_{\vartriangle}\{a\}\bigr)\geq h\ \ \forall
a\in A\bigr\}$ is described in~(\ref{eq:18}). Further, for any
$h\geq\nu(r\cdot k)$, we by~(\ref{eq:15}) have
$\mathfrak{F}\bigl(\boldsymbol{\mathsf{b}}_{\nu(r\cdot
k)-1}(A)\bigr)\supseteq
\mathfrak{F}\bigl(\boldsymbol{\mathsf{b}}_{h-1}(A)\bigr)$,
so~(\ref{eq:21}) reduces to
\begin{equation*}
\mathbf{I}_r(P,A;\omega)=\bigcup_{1\leq k\leq \omega(P)}
\bigl\{b\in P:\ \omega(b)=k,\
\omega\bigl(\{b\}\wedge_{\vartriangle}\{a\}\bigr)\geq \nu(r\cdot
k)\ \ \forall a\in A\bigr\}\ ,
\end{equation*}
and we come to the following conclusion.

\begin{prop}\label{prop:4} Let $A$ be a nontrivial antichain in $P$.
\begin{itemize}
\item[\rm(i)]
For any map $\omega$, it holds
\begin{equation*}
\begin{split}
\mathbf{I}_r(P,A;\omega)&= \bigcup_{1\leq k\leq \omega(P)}
\biggl(\Bigl(
\mathfrak{F}\bigl(\boldsymbol{\mathsf{b}}_{k-1}(\hat{1}_P)\bigr)-\mathfrak{F}
\bigl(\boldsymbol{\mathsf{b}}_k(\hat{1}_P)\bigr) \Bigr)\cap
\mathfrak{F}\bigl(\boldsymbol{\mathsf{b}}_{\nu(r\cdot
k)-1}(A)\bigr)\biggr)\\&= \bigcup_{1\leq k\leq \omega(P)}\Bigl(\
\mathfrak{F}\bigl(\boldsymbol{\mathsf{b}}_{k-1}(\hat{1}_P)\wedge_{\triangledown}
\boldsymbol{\mathsf{b}}_{\nu(r\cdot k)-1}(A)\bigr)-
\mathfrak{F}\bigl(\boldsymbol{\mathsf{b}}_{k}(\hat{1}_P)\bigr)\
\Bigr)\ .
\end{split}
\end{equation*}
\item[\rm(ii)] If $P$ is graded, then
\begin{equation}
\label{eq:22} \mathbf{I}_r(P,A;\rho)= \bigcup_{k\in\lbrack
1,\rho(P)\rbrack:\ \nu(r\cdot k)\leq\min_{a\in A}\rho(a)} \Bigl(
P^{(k)}\cap\mathfrak{F}\bigl(\boldsymbol{\mathsf{b}}_{\nu(r\cdot
k)-1}(A)\bigr)\Bigr)\ .
\end{equation}
\end{itemize}
\end{prop}

To find the cardinality of subposet~(\ref{eq:22}), we can use the
combinatorial inclusion-exclusion principle (see, e.g.,~Chapter~IV
of Aigner, 1979, and Chapter~2 of Stanley, 1997). Under the
hypothesis of Proposition~\ref{prop:4}(ii), we have
\begin{align*}
|\mathbf{I}_r(P,A;\rho)|&=\sum_{k\in\lbrack 1,\rho(P)\rbrack:\
\nu(r\cdot k)\leq\min_{a\in A}\rho(a)}\\&\phantom{=}
\sum_{C\subseteq A:\ |C|>0 }(-1)^{|C|-1}\cdot\left|
P^{(k)}\cap\mathfrak{F}\left(\mathfrak{I}(C)\cap P^{(\nu(r\cdot
k))}\right)\right|\\&= \sum_{k\in\lbrack 1,\rho(P)\rbrack:\
\nu(r\cdot k)\leq\min_{a\in
A}\rho(a)}\\&\phantom{=}\sum_{E\subseteq P^{(\nu(r\cdot
k))}\cap\mathfrak{I}(A):\ |E|>0} \left(\sum_{C\subseteq A:\
E\subseteq\mathfrak{I}(C)}(-1)^{|C|-1}\right)\cdot \left|
P^{(k)}\cap\mathfrak{F}(E)\right| \ .
\end{align*}

For the remainder of the present section, let $A$ be a nontrivial
antichain in a graded lattice $P$ of rank $n$, with the property:
each interval of length $k$ in $P$ contains the same number
$\mathrm{B}(k)$ of maximal chains; in other words, we suppose $P$
to be a principal order ideal of some {\em binomial poset},
see~Section~3.15 of Stanley, 1997. The function $\mathrm{B}(k)$ is
called the {\em binomial function\/} of $P$; it holds
$\mathrm{B}(0)=\mathrm{B}(1)=1$. The number of elements of rank
$i$ in any interval of length $j$ is denoted by
$\left[\begin{smallmatrix}j\\ i\end{smallmatrix}\right]$; it holds
$\left[\begin{smallmatrix}j\\
i\end{smallmatrix}\right]=\tfrac{\mathrm{B}(j)}{\mathrm{B}(i)\cdot\mathrm{B}(j-i)}$.
If $P$ is $\mathbb{B}(n)$ or $\mathbb{V}_q(n)$, then
$\left[\begin{smallmatrix}j\\
i\end{smallmatrix}\right]=\tbinom{j}{i}$ or
$\left[\begin{smallmatrix}j\\
i\end{smallmatrix}\right]=\tbinom{j}{i}_{\! q}$, respectively.

Given $k\in\lbrack 1,n\rbrack$ such that $\nu(r\cdot
k)\leq\min_{a\in A}\rho(a)$, we have
\begin{equation}\label{eq:23}
\begin{split}
|\mathbf{I}_{r,k}(P,A;\rho)|& =\sum_{C\subseteq A:\
|C|>0}(-1)^{|C|-1}\\&\phantom{=}\cdot\sum_{E\subseteq
P^{(\nu(r\cdot k))}\cap\mathfrak{I}(C):\
|E|>0}(-1)^{|E|-1}\cdot\begin{bmatrix}n-\rho\left(\bigvee_{e\in
E}e\right)\\n-k\end{bmatrix}\\&= \sum_{E\subseteq P^{(\nu(r\cdot
k))}\cap\mathfrak{I}(A):\
|E|>0}(-1)^{|E|}\\&\phantom{=}\cdot\left(\sum_{C\subseteq A:\
E\subseteq\mathfrak{I}(C)}(-1)^{|C|}\right)\cdot \begin{bmatrix}
n-\rho\left(\bigvee_{e\in E}e\right)\\n-k\end{bmatrix} \ .
\end{split}
\end{equation}
Indeed, for example, the sum
\begin{equation}\label{eq:24}
\sum_{E\subseteq P^{(\nu(r\cdot k))}\cap\mathfrak{I}(C):\
|E|>0}(-1)^{|E|-1}\cdot\begin{bmatrix}n-\rho\left(\bigvee_{e\in
E}e\right)\\n-k\end{bmatrix}
\end{equation}
counts the number of elements of the layer $P^{(k)}$ comparable
with, at least, one element of the antichain $P^{(\nu(r\cdot
k))}\cap\mathfrak{I}(C)$.

To refine expression~(\ref{eq:23}) with the help of the technique
of the M\"{o}bius function, consider some auxiliary lattices which
can be associated to the antichain $A$. The first one, denoted by
$\mathcal{C}_{r,k}(P,A)$, is the lattice consisting of all sets
from the family $\{P^{(\nu(r\cdot k))}\cap\mathfrak{I}(C):\
C\subseteq A\}$ ordered by inclusion. The greatest element of
$\mathcal{C}_{r,k}(P,A)$ is the set $P^{(\nu(r\cdot
k))}\cap\mathfrak{I}(A)$. The least element of
$\mathcal{C}_{r,k}(P,A)$, denoted by $\hat{0}$, is the empty
subset of $P^{(\nu(r\cdot k))}$. The remaining lattices, denoted
by $\mathcal{E}_{r,k}(P,X)$, where $X$ are nonempty subsets of
$P^{(\nu(r\cdot k))}\cap\mathfrak{I}(A)$, are defined in the
following way. Given an antichain $X\subseteq P^{(\nu(r\cdot
k))}\cap\mathfrak{I}(A)$, the poset $\mathcal{E}_{r,k}(P,X)$ is
the sub-join-semilattice of the lattice $P$ generated by $X$ and
augmented with a new least element, denoted by $\hat{0}$ (it is
regarded as the empty subset of $P$). The greatest element of
$\mathcal{E}_{r,k}(P,X)$ is the join $\bigvee_{x\in X}x$ in $P$.
We have
\begin{equation}\label{eq:25}
\begin{split}
|\mathbf{I}_{r,k}(P,A;\rho)| &=\sum_{X\in\mathcal{C}_{r,k}(P,A):\
\hat{0}<X}\mu_{\mathcal{C}_{r,k}(P,A)}(\hat{0},X)\\&\phantom{=}\cdot
\sum_{z\in\mathcal{E}_{r,k}(P,X):\ \hat{0}<z,\ \rho(z)\leq k
}\mu_{\mathcal{E}_{r,k}(P,X)}(\hat{0},z)\cdot
\begin{bmatrix}n-\rho(z)\\n-k\end{bmatrix} \ ,
\end{split}
\end{equation}
where $\rho(\cdot)$ means the rank in $P$, and where, for example,
the sum
\begin{equation*}
-\sum_{z\in\mathcal{E}_{r,k}(P,X):\ \hat{0}<z,\ \rho(z)\leq
k}\mu_{\mathcal{E}_{r,k}(P,X)}(\hat{0},z)\cdot
\begin{bmatrix}n-\rho(z)\\n-k\end{bmatrix}
\end{equation*}
is equivalent to sum~(\ref{eq:24}) under $X=P^{(\nu(r\cdot
k))}\cap\mathfrak{I}(C)$.

If $P$ is $\mathbb{B}(n)$ then, in view of
Remark~{\rm\ref{rem:4}}, formulas~{\rm(\ref{eq:23})}
and~{\rm(\ref{eq:25})} give, for a nontrivial clutter, the number
of all its $r$-committees of cardinality $k$.

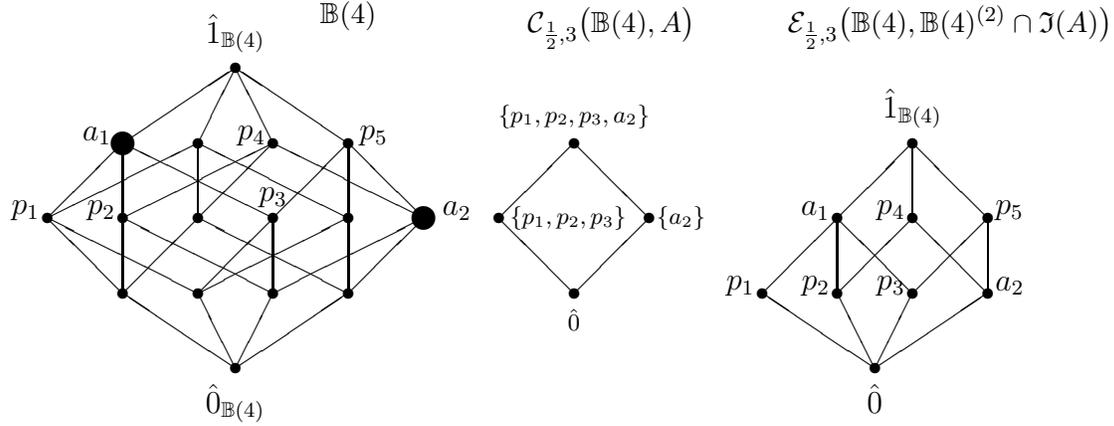
\begin{figure}[h]
\begin{picture}(11,6)(2,0)

\put(0.5,3){\circle * {0.15}} \put(1.5,2){\circle * {0.15}}
\put(1.5,3){\circle
* {0.15}} \put(3,1){\circle * {0.15}} \put(2.5,2){\circle * {0.15}}
\put(2.5,3){\circle * {0.15}} \put(3,5){\circle * {0.15}}
\put(3.5,2){\circle * {0.15}} \put(3.5,3){\circle * {0.15}}
\put(5.5,3){\circle * {0.3}}

\put(4.5,2){\circle * {0.15}} \put(4.5,3){\circle * {0.15}}

\put(1.5,4){\circle * {0.3}} \put(2.5,4){\circle * {0.15}}
\put(3.5,4){\circle
* {0.15}} \put(4.5,4){\circle * {0.15}}

\put(3,1){\line (-3,2){1.5}} \put(3,1){\line (-1,2){0.5}}
\put(3,1){\line (1,2){0.5}} \put(3,1){\line (3,2){1.5}}
\put(1.5,2){\line (-1,1){1}} \put(1.5,2){\line (0,1){1}}
\put(1.5,2){\line (1,1){1}}

\put(0.5,3){\line (1,1){1}} \put(0.5,3){\line (2,1){2}}

\put(1.5,4){\line (3,2){1.5}}

\put(2.5,2){\line (-2,1){2}} \put(2.5,2){\line (1,1){1}}
\put(2.5,2){\line (2,1){2}}

\put(3.5,2){\line (-2,1){2}} \put(3.5,2){\line (0,1){1}}
\put(3.5,2){\line (2,1){2}}

\put(4.5,2){\line (-2,1){2}} \put(4.5,2){\line (0,1){1}}
\put(4.5,2){\line (1,1){1}}

\put(1.5,3){\line (0,1){1}} \put(1.5,3){\line (2,1){2}}

\put(2.5,3){\line (0,1){1}} \put(2.5,3){\line (1,1){1}}

\put(3.5,3){\line (-2,1){2}} \put(3.5,3){\line (1,1){1}}

\put(4.5,3){\line (-2,1){2}} \put(4.5,3){\line (0,1){1}}

\put(5.5,3){\line (-2,1){2}} \put(5.5,3){\line (-1,1){1}}

\put(2.5,4){\line (1,2){0.5}}

\put(3.5,4){\line (-1,2){0.5}}

\put(4.5,4){\line (-3,2){1.5}}

\put(3,0.8){\makebox(0,0)[t]{$\hat{0}_{\mathbb{B}(4)}$}}
\put(3,5.2){\makebox(0,0)[b]{$\hat{1}_{\mathbb{B}(4)}$}}

\put(1.35,4.1){\makebox(0,0)[r]{$a_1$}}
\put(3.4,4.1){\makebox(0,0)[r]{$p_4$}}
\put(4.65,4.1){\makebox(0,0)[l]{$p_5$}}
\put(5.75,3.1){\makebox(0,0)[l]{$a_2$}}

\put(0.4,3.1){\makebox(0,0)[r]{$p_1$}}
\put(1.4,3.1){\makebox(0,0)[r]{$p_2$}}
\put(3.5,3.15){\makebox(0,0)[b]{$p_3$}}

\put(4.5,5.5){\makebox(0,0)[b]{\small $\mathbb{B}(4)$}}

\put(6.5,3){\line (1,1){1}} \put(7.5,4){\line (1,-1){1}}
\put(7.5,2){\line (-1,1){1}} \put(7.5,2){\line (1,1){1}}

\put(6.5,3){\circle * {0.15}} \put(8.5,3){\circle *{0.15}}
\put(7.5,4){\circle * {0.15}} \put(7.5,2){\circle * {0.15}}

\put(6.65,3){\makebox(0,0)[l]{\scriptsize $\{p_1,p_2,p_3\}$}}
\put(8.6,3){\makebox(0,0)[l]{\scriptsize $\{a_2\}$}}
\put(7.5,4.2){\makebox(0,0)[b]{\scriptsize $\{p_1,p_2,p_3,a_2\}$}}
\put(7.5,1.8){\makebox(0,0)[t]{\scriptsize $\hat{0}$}}

\put(8,5.25){\makebox(0,0)[b]{\small $\mathcal{C}_{\frac{1}{2},3}
\bigl(\mathbb{B}(4),A\bigr)$}}


\put(10,2){\circle * {0.15}} \put(11.5,1){\circle * {0.15}}
\put(11,2){\circle
* {0.15}} \put(11,3){\circle * {0.15}} \put(12,2){\circle * {0.15}}
\put(12,3){\circle * {0.15}} \put(13,2){\circle * {0.15}}
\put(13,3){\circle
* {0.15}} \put(12,4){\circle * {0.15}}

\put(11.5,1){\line (-3,2){1.5}} \put(11.5,1){\line (-1,2){0.5}}
\put(11.5,1){\line (1,2){0.5}} \put(11.5,1){\line (3,2){1.5}}

\put(10,2){\line (1,1){2}} \put(11,2){\line (0,1){1}}
\put(11,2){\line (1,1){1}} \put(12,2){\line (-1,1){1}}
\put(12,2){\line (1,1){1}}\put(13,2){\line (-1,1){1}}
\put(13,2){\line (0,1){1}} \put(12,3){\line (0,1){1}}
\put(13,3){\line (-1,1){1}}

\put(11.5,0.8){\makebox(0,0)[t]{$\hat{0}$}}

\put(9.9,2.1){\makebox(0,0)[r]{$p_1$}}
\put(10.9,2.1){\makebox(0,0)[r]{$p_2$}}
\put(11.9,2.1){\makebox(0,0)[r]{$p_3$}}
\put(13.1,2.1){\makebox(0,0)[l]{$a_2$}}
\put(10.9,3.1){\makebox(0,0)[r]{$a_1$}}
\put(11.9,3.1){\makebox(0,0)[r]{$p_4$}}
\put(13.1,3.1){\makebox(0,0)[l]{$p_5$}}

\put(12,4.2){\makebox(0,0)[b]{$\hat{1}_{\mathbb{B}(4)}$}}

\put(12.5,5.25){\makebox(0,0)[b]{\small
$\mathcal{E}_{\frac{1}{2},3}\!
\left(\mathbb{B}(4),\mathbb{B}(4)^{(2)}\cap\mathfrak{I}(A)\right)$}}

\end{picture}
\caption{An antichain in the Boolean lattice, and auxiliary
lattices involved in enumeration of relatively blocking elements}
\label{illustration}
\end{figure}

{\small
\begin{exmpl} Figure~{\rm\ref{illustration}} depicts the Hasse diagram of
a Boolean lattice of rank four, its antichain $A:=\{a_1,a_2\}$,
and lattices
$\mathcal{C}:=\mathcal{C}_{\frac{1}{2},3}(\mathbb{B}(4),A)$ and
$\mathcal{E}:=\mathcal{E}_{\frac{1}{2},3}
(\mathbb{B}(4),\mathbb{B}(4)^{(2)}\cap\mathfrak{I}(A))$. To
compute the number of elements in
$\mathbf{I}_{\frac{1}{2},3}(\mathbb{B}(4),A;\rho)$, note that
$\mu_{\mathcal{C}}(\hat{0},\{p_1,p_2,p_3\})=
\mu_{\mathcal{C}}(\hat{0},\{a_2\})=-1$ and
$\mu_{\mathcal{C}}(\hat{0},\{p_1,p_2,p_3,a_2\})=1$. Further, we
have
$\mu_{\mathcal{E}}(\hat{0},p_1)=\mu_{\mathcal{E}}(\hat{0},p_2)=
\mu_{\mathcal{E}}(\hat{0},p_3)=\mu_{\mathcal{E}}(\hat{0},a_2)=-1$,
$\mu_{\mathcal{E}}(\hat{0},a_1)=2$,
$\mu_{\mathcal{E}}(\hat{0},p_4)=\mu_{\mathcal{E}}(\hat{0},p_5)=1$
and $\mu_{\mathcal{E}}(\hat{0},\hat{1}_{\mathbb{B}(4)})=-1$.

By means of~{\rm(\ref{eq:25})}, we obtain
$|\mathbf{I}_{\frac{1}{2},3}(\mathbb{B}(4),A;\rho)|=|\{p_4,p_5\}|=2$.
\end{exmpl}
}

Proposition~\ref{prop:4}(ii) provides us with a general
description of the subposets of relatively $r$-blocking elements
in graded posets. The aim of
Section~\ref{SectionOnApplicationOfFareySubsequences} of the
present paper is to explore the structure of the above-mentioned
subposets in detail; with the help of Theorem~\ref{thm:1} we will
exclude from consideration some layers of graded posets that
certainly contain no relatively $r$-blocking elements.

\section{Principal order ideals and Farey subsequences}
\label{SectionOnFareySubsequences}

Let $\{a\}$ be a one-element antichain in $P$. Define the sequence
of irreducible fractions
\begin{multline*}
\mathcal{F}(P,a;\omega):=\left\{\frac{0}{1},\frac{1}{1}\right\}\\
\cup\left(\frac{\omega(\{b\}\wedge_{\vartriangle}\{a\})}
{\gcd\bigl(\omega(\{b\}\wedge_{\vartriangle}\{a\}),\omega(b)\bigr)}\Bigg{/}
\frac{\omega(b)}
{\gcd\bigl(\omega(\{b\}\wedge_{\vartriangle}\{a\}),\omega(b)\bigr)}:\
b\in P-\{\hat{0}_P\}\right)
\end{multline*}
arranged in ascending order.

Recall that the {\em Farey sequence\/} $\mathcal{F}_n$ of order
$n\in\mathbb{P}$ is defined to be the ascending sequence of all
irreducible fractions between $0$ and $1$ whose denominators do
not exceed $n$, see, e.g.,~Chapter~27 of Buchstab, 1967, Chapter~4
of Graham et al., 1994, Chapter~III of Hardy and Wright, 1979, and
Lagarias and Tresser, 1995. Thus, $\mathcal{F}(P,a;\omega)$ is a
subsequence of the Farey sequence of order $\omega(P)$.

We always index the fractions from $\mathcal{F}(P,a;\omega)$
starting with zero:
$\mathcal{F}(P,a;\omega)=\bigl(f_0:=\tfrac{0}{1}<f_1<f_2\cdots
<f_{|\mathcal{F}(P,a;\omega)|-1}:=\tfrac{1}{1}\bigr)$.

In the present paper, we do not deal with the more general
ascending Farey subsequences $\bigcup_{a\in
A}\mathcal{F}(P,a;\omega)$ and $\bigcap_{a\in
A}\mathcal{F}(P,a;\omega)$ associated to nonempty antichains $A$
in $P$; such sequences can also be of interest.

Order-preserving maps $P\to\mathbb{P}$ and
$P\to\mathbb{P}^{\ast}$, where $\mathbb{P}^{\ast}$ are positive
integers ordered by divisibility, are discussed, e.g., in~Smith,
1967, 1969, 1970/1971.

See~P\v{a}tra\c{s}cu and P\v{a}tra\c{s}cu, 2004, on algorithmic
aspects of the Farey sequences.

\section{Farey subsequences in Boolean context}
\label{SectionOnFareySubsequencesInBooleanContext}

In this section we deal almost exclusively with the Boolean
lattice $\mathbb{B}(n)$. Let $a$ be an arbitrary element of
$\mathbb{B}(n)$, of rank $m:=\rho(a)$. Consider the Farey
subsequence $\mathcal{F}(\mathbb{B}(n),a;\rho)$ associated to the
principal order ideal $\mathfrak{I}(a)$ of $\mathbb{B}(n)$. The
sequences $\mathcal{F}(\mathbb{B}(n),a;\rho)$ are the same, for
all elements $a$ of rank $m$ in $\mathbb{B}(n)$, and we write
$\mathcal{F}(\mathbb{B}(n),m;\rho)$ instead of
$\mathcal{F}(\mathbb{B}(n),a;\rho)$. For any element
$b\in\mathbb{B}(n)-\{\hat{0}_{\mathbb{B}(n)}\}$, we have
$\rho(a)+\rho(b)-\rho(a\wedge b)=m+\rho(b)-\rho(a\wedge b)\leq n$;
moreover, $0\leq\rho(a\wedge b)\leq\rho(a)=:m$, so we are
interested in the ascending Farey subsequence
\begin{equation}\label{eq:26}
\mathcal{F}(\mathbb{B}(n),m;\rho)=\left\{\tfrac{0}{1},\tfrac{1}{1}\right\}\cup\left(\tfrac{h}{k}\in\mathcal{F}_n:\
h\leq m,\ k-h\leq n-m\right)\ .
\end{equation}

{\small
\begin{exmpl}
\label{exmpl:1}
\begin{equation*}
\begin{split}
\mathcal{F}_6&= \left(
\tfrac{0}{1}<\tfrac{1}{6}<\tfrac{1}{5}<\tfrac{1}{4}<\tfrac{1}{3}<\tfrac{2}{5}<
\tfrac{1}{2}<\tfrac{3}{5}<\tfrac{2}{3}<\tfrac{3}{4}<\tfrac{4}{5}<\tfrac{5}{6}<
\tfrac{1}{1} \right)\ ,\\ \mathcal{F}(\mathbb{B}(6),6;\rho)&=
\left( \tfrac{0}{1}<\tfrac{1}{1} \right)\ ,\\
\left(\tfrac{h}{k}\in\mathcal{F}_6:\ h\leq 5\right) &=
\mathcal{F}_6\ ,\\ \mathcal{F}(\mathbb{B}(6),5;\rho)&=\left(
\tfrac{0}{1}<
\tfrac{1}{2}<\tfrac{2}{3}<\tfrac{3}{4}<\tfrac{4}{5}<\tfrac{5}{6}<
\tfrac{1}{1} \right)\ ,\\ \left(\tfrac{h}{k}\in\mathcal{F}_6:\
h\leq 4\right) &= \left(
\tfrac{0}{1}<\tfrac{1}{6}<\tfrac{1}{5}<\tfrac{1}{4}<\tfrac{1}{3}<\tfrac{2}{5}<
\tfrac{1}{2}<\tfrac{3}{5}<\tfrac{2}{3}<\tfrac{3}{4}<\tfrac{4}{5}<\tfrac{1}{1}
\right)\ ,\\ \mathcal{F}(\mathbb{B}(6),4;\rho)&=\left(
\tfrac{0}{1}<\tfrac{1}{3}<
\tfrac{1}{2}<\tfrac{3}{5}<\tfrac{2}{3}<\tfrac{3}{4}<\tfrac{4}{5}<\tfrac{1}{1}
\right)\ ,\\ \left(\tfrac{h}{k}\in\mathcal{F}_6:\ h\leq 3\right)
&= \left(
\tfrac{0}{1}<\tfrac{1}{6}<\tfrac{1}{5}<\tfrac{1}{4}<\tfrac{1}{3}<\tfrac{2}{5}<
\tfrac{1}{2}<\tfrac{3}{5}<\tfrac{2}{3}<\tfrac{3}{4}< \tfrac{1}{1}
\right)\ ,
\\ \mathcal{F}(\mathbb{B}(6),3;\rho)&= \left(
\tfrac{0}{1}<\tfrac{1}{4}<\tfrac{1}{3}<\tfrac{2}{5}<
\tfrac{1}{2}<\tfrac{3}{5}<\tfrac{2}{3}<\tfrac{3}{4}< \tfrac{1}{1}
\right)\ ,\\ \left(\tfrac{h}{k}\in\mathcal{F}_6:\ h\leq 2\right)
&= \left(
\tfrac{0}{1}<\tfrac{1}{6}<\tfrac{1}{5}<\tfrac{1}{4}<\tfrac{1}{3}<\tfrac{2}{5}<
\tfrac{1}{2}<\tfrac{2}{3}<\tfrac{1}{1} \right)\ ,\\
\mathcal{F}(\mathbb{B}(6),2;\rho)&= \left(
\tfrac{0}{1}<\tfrac{1}{5}<\tfrac{1}{4}<\tfrac{1}{3}<\tfrac{2}{5}<
\tfrac{1}{2}<\tfrac{2}{3}<\tfrac{1}{1} \right)\ ,\\
\left(\tfrac{h}{k}\in\mathcal{F}_6:\ h\leq 1\right) &= \left(
\tfrac{0}{1}<\tfrac{1}{6}<\tfrac{1}{5}<\tfrac{1}{4}<\tfrac{1}{3}<
\tfrac{1}{2}< \tfrac{1}{1} \right)\ ,\\
\mathcal{F}(\mathbb{B}(6),1;\rho)&=\left(\tfrac{h}{k}\in\mathcal{F}_6:\
h\leq 1\right)\ ,\\ \left(\tfrac{h}{k}\in\mathcal{F}_6:\ h=
0\right) &= \left( \tfrac{0}{1}\right)\ ,\\
\mathcal{F}(\mathbb{B}(6),0;\rho)&= \left(
\tfrac{0}{1}<\tfrac{1}{1} \right)\ .
\end{split}
\end{equation*}
\end{exmpl}
}

{\small
\begin{rem} \label{rem:5}
In the sequence $\mathcal{F}(\mathbb{B}(n),m;\rho)$ such that
$0<m<n$, we have
\begin{equation*}
f_0=\tfrac{0}{1},\ \ f_1=\tfrac{1}{n-m+1},\ \
f_{|\mathcal{F}(\mathbb{B}(n),m;\rho)|-2}=\tfrac{m}{m+1},\ \
f_{|\mathcal{F}(\mathbb{B}(n),m;\rho)|-1}=\tfrac{1}{1}.
\end{equation*}
\end{rem}
}

Let $n\in\mathbb{P}$, and let $S$ be a subset of $\lbrack
1,n\rbrack$. We denote by $\phi(n;S)$ the number of elements from
$S$ that are relatively prime to $n$:
\begin{equation}\label{eq:28}
\phi(n;S):=|\{s\in S:\ n\bot s\}|\ ;
\end{equation}
thus, $\phi\bigl(n;\lbrack 1,n\rbrack\bigr)$ is the {\em Euler
function}.

Given a positive integer $i$ such that $i\leq n$, we have
$\phi\bigl(n;\lbrack 1,i\rbrack\bigr)=\sum_{d\in\lbrack
1,i\rbrack:\
d|n}\overline{\mu}(d)\cdot\left\lfloor\frac{i}{d}\right\rfloor$,
where $\overline{\mu}(\cdot)$ stands for the number-theoretic {\em
M\"{o}bius function}: $\overline{\mu}(1):=1$; if $p^2|d$, for some
prime $p$, then $\overline{\mu}(d):=0$; if $d=p_1p_2\cdots p_s$,
for distinct primes $p_1,p_2,\ldots,p_s$, then
$\overline{\mu}(d):=(-1)^s$. Thus, given a nonempty subset
$\lbrack i'+1,i''\rbrack\subseteq\lbrack 1,n\rbrack$, we have
$\phi\bigl(n;\lbrack i'+1,i''\rbrack\bigr)=\sum_{d\in\lbrack
1,i''\rbrack:\ d|n}\overline{\mu}(d)\cdot \left(
\left\lfloor\frac{i''}{d}\right\rfloor -
\left\lfloor\frac{i'}{d}\right\rfloor \right)$.

\begin{prop}\label{prop:5}
\begin{itemize}
\item[\rm(i)]
If
$f_t\in\mathcal{F}(\mathbb{B}(n),m;\rho)-\left\{\tfrac{1}{1}\right\}$,
where $0<m<n$, then
\begin{equation*}
t=\sum_{j\in\lbrack 1,n\rbrack}\phi \Bigl(
j;\Bigl\lbrack\max\bigl\{1,j+\min\{m,\lfloor j\cdot
f_t\rfloor\}-n\bigr\},\ \min\{m,\lfloor j\cdot
f_t\rfloor\}\Bigr\rbrack \Bigr)\ .
\end{equation*}
\item[\rm(ii)] The cardinality of the
sequence $\mathcal{F}(\mathbb{B}(n),m;\rho)$, where $0<m<n$,
equals
\begin{multline*}
1+\sum_{j\in\lbrack 1,n\rbrack}\phi\bigl(j;\bigl\lbrack
1,\min\{m,j\}\bigr\rbrack\bigr)-\sum_{j\in\left\lbrack\left\lceil
n/2\right\rceil+1,m\right\rbrack}\phi\bigl(j;\lbrack 1,2\cdot
j-n-1 \rbrack\bigr)\\-\sum_{j\in\lbrack
n-m+2,n\rbrack}\phi\bigl(j;\lbrack 1,j+m-n-1 \rbrack\bigr)\ .
\end{multline*}
\end{itemize}
\end{prop}

\begin{proof} To prove assertion~(i),
replace $f_t$ with $\tfrac{j\cdot f_t}{j}$, for every $j\in\lbrack
1,n\rbrack$. According to description~(\ref{eq:26}), $t$ equals
\begin{multline*}
\sum_{j\in\lbrack 1,n\rbrack} \bigl| \bigl\{\ i\in\lbrack
1,j\rbrack:\\ i\bot j,\ \max\bigl\{1,j+ \min\{m,\lfloor j\cdot
f_t\rfloor\}-n\bigr\} \leq i\leq\min\{m,\lfloor j\cdot
f_t\rfloor\}\ \bigr\}\bigr| \ ,
\end{multline*}
from where the assertion follows, with respect to
description~(\ref{eq:28}) of the function $\phi(\cdot;\cdot)$.
Assertion~(i) implies assertion~(ii), due to our convention that
$\tfrac{1}{1}$ is the terminal fraction in the sequence
$\mathcal{F}(\mathbb{B}(n),m;\rho)$, with the index
$|\mathcal{F}(\mathbb{B}(n),m;\rho)|-1$. Indeed, we have
\begin{multline*} |\mathcal{F}(\mathbb{B}(n),m;\rho)|-1=
\sum_{j\in\lbrack 1,n\rbrack}
\phi\left(j;\bigl\lbrack\max\{1,j+\min\{m,j\}-n\},\
\min\{m,j\}\bigr\rbrack\right)
\\
=\sum_{j\in\lbrack
1,m\rbrack}\phi\bigl(j;\bigl\lbrack\max\{1,2\cdot j-n\},\
j\bigr\rbrack\bigr)+ \sum_{j\in\lbrack m+1,n\rbrack}
\phi\bigl(j;\bigl\lbrack\max\{1,j+m-n\},\ m\bigr\rbrack\bigr)\ ,
\end{multline*}
and assertion~(ii) follows.
\end{proof}

The sum $1+\sum_{j\in\lbrack 1,n\rbrack}\phi\bigl(j;\bigl\lbrack
1,\min\{m,j\}\bigr\rbrack\bigr)$ appearing in
Proposition~{\rm\ref{prop:5}(ii)} counts the number of fractions
in the sequence $(\tfrac{h}{k}\in\mathcal{F}_n:\ h\leq m)$, see
Remark~{\rm\ref{rem:6}(ii)(b)} below.

Description~(\ref{eq:26}) of Farey subsequences leads up to the
following observation.

\begin{prop}
The map
\[
\mathcal{F}(\mathbb{B}(n),m;\rho)\to\mathcal{F}(\mathbb{B}(n),n-m;\rho)\
,\ \ \tfrac{h}{k}\mapsto\tfrac{k-h}{k}
\]
is order-reversing and bijective, for any\/ $m$, $0\leq m\leq n$.
\end{prop}

We now explore the properties of Farey subsequences~(\ref{eq:26}).

\begin{prop}
\label{prop:6} Let
$\frac{h}{k}\in\mathcal{F}(\mathbb{B}(n),m;\rho)$, where $0<m<n$;
suppose that $\frac{0}{1}<\frac{h}{k}<\frac{1}{1}$.
\begin{itemize}
\item[\rm(i)]
Let $x_0$ be the integer such that $kx_0\equiv -1\pmod{h}$ and
$m-h+1\leq x_0\leq m$. Define integers $y_0$ and $t^{\ast}$ by
$y_0:=\frac{kx_0+1}{h}$ and $t^\ast:=
\left\lfloor\min\{\frac{m-x_0}{h},\frac{n-y_0}{k},
\frac{n-m+x_0-y_0}{k-h}\}\right\rfloor$.

The fraction $\frac{x_0+t^\ast h}{y_0+t^\ast k}$ precedes the
fraction $\frac{h}{k}$ in $\mathcal{F}(\mathbb{B}(n),m;\rho)$.

\item[\rm(ii)]
Let $x_0$ be the integer such that $kx_0\equiv 1\pmod{h}$ and
$m-h+1\leq x_0\leq m$. Define integers $y_0$ and $t^{\ast}$ by
$y_0:=\frac{kx_0-1}{h}$ and $t^{\ast}:=
\left\lfloor\min\{\frac{m-x_0}{h},\frac{n-y_0}{k},
\frac{n-m+x_0-y_0}{k-h}\}\right\rfloor$.

The fraction $\frac{x_0+t^{\ast} h}{y_0+t^{\ast} k}$ succeeds the
fraction $\frac{h}{k}$ in $\mathcal{F}(\mathbb{B}(n),m;\rho)$.
\end{itemize}
\end{prop}

\begin{sketch}
We sketch the proof of assertion~(i).

Since the pair $(x_0,y_0)$ is a solution to the equation
$-kx+hy=1$, the pair $(x_0+th,y_0+tk)$ is a solution as well, for
any integer $t$. Considering the system of inequalities $0\leq
x_0+th\leq m$, $1\leq y_0+tk\leq n$, $1\leq y_0+tk-(x_0+th)\leq
n-m$, where $t$ is an integer variable, we can turn to the
solution-equivalent system
\begin{equation}
\label{eq:29}
\begin{cases}
-\frac{x_0}{h}\leq t\leq\frac{m-x_0}{h}\ ,\\ \frac{-y_0+1}{k}\leq
t\leq\frac{n-y_0}{k}\ ,\\ \frac{x_0-y_0+1}{k-h}\leq
t\leq\frac{n-m+x_0-y_0}{k-h}\ .
\end{cases}
\end{equation}
Note that
$\max\left\{-\frac{x_0}{h},\frac{-y_0+1}{k},\frac{x_0-y_0+1}{k-h}
\right\}=\frac{x_0-y_0+1}{k-h}$, therefore system~(\ref{eq:29}) is
solution-equivalent to the inequality
\begin{equation}\label{eq:30}
\frac{x_0-y_0+1}{k-h}\leq t\leq\min\left\{\frac{m-x_0}{h},
\frac{n-y_0}{k},\frac{n-m+x_0-y_0}{k-h}\right\}\ .
\end{equation}
Inequality~(\ref{eq:30}) has at least one integer solution, namely
$t=\left\lceil\frac{x_0-y_0+1}{k-h}\right\rceil$. Another
observation is that, for any integer solutions $t'$ and $t''$
to~(\ref{eq:30}) such that $t'\leq t''$, we have
$\tfrac{0}{1}\leq\tfrac{x_0+t'h}{y_0+t'k}\leq\tfrac{x_0+t''h}{y_0+t''k}<\tfrac{h}{k}$.
The proof of assertion~(i) is completed by checking that there is
no fraction $\tfrac{i}{j}\in\mathcal{F}(\mathbb{B}(n),m;\rho)$
such that
$\tfrac{x_0+t^{\ast}h}{y_0+t^{\ast}k}<\tfrac{i}{j}<\tfrac{h}{k}$;
thus, the fraction $\tfrac{x_0+t^{\ast}h}{y_0+t^{\ast}k}$ does
precede the fraction $\tfrac{h}{k}$ in
$\mathcal{F}(\mathbb{B}(n),m;\rho)$.

Assertion~(ii) can be proved in an analogous way.
\end{sketch}

{\small
\begin{rem}\label{rem:3}
If $\tfrac{1}{k}\in\mathcal{F}(\mathbb{B}(n),m;\rho)$, for some
$k>1$, then Proposition~{\rm\ref{prop:6}} implies that the
fraction
\begin{equation*}
\frac{m+\min\left\{0,\left\lfloor\tfrac{n-km-1}{k-1}\right\rfloor\right\}}
{k\cdot\left(m+
\min\left\{0,\left\lfloor\tfrac{n-km-1}{k-1}\right\rfloor\right\}\right)+1}
\end{equation*}
precedes $\tfrac{1}{k}$, and the fraction
\begin{equation*}
\frac{m+\min\left\{0,\left\lfloor\tfrac{n-km+1}{k-1}\right\rfloor\right\}}
{k\cdot\left(m+
\min\left\{0,\left\lfloor\tfrac{n-km+1}{k-1}\right\rfloor\right\}\right)-1}
\end{equation*}
succeeds $\tfrac{1}{k}$ in $\mathcal{F}(\mathbb{B}(n),m;\rho)$.
\end{rem}
}

\begin{prop}
\label{prop:7}
\begin{itemize}
\item[\rm(i)]
If $\frac{h_j}{k_j}<\frac{h_{j+1}}{k_{j+1}}$ are two successive
fractions of $\mathcal{F}(\mathbb{B}(n),m;\rho)$, where $0\leq
m\leq n$, then
\begin{equation}\label{eq:31}
k_j h_{j+1}-h_j k_{j+1}=1\ .
\end{equation}

\item[\rm(ii)]
If
$\frac{h_j}{k_j}<\frac{h_{j+1}}{k_{j+1}}<\frac{h_{j+2}}{k_{j+2}}$
are three successive fractions of
$\mathcal{F}(\mathbb{B}(n),m;\rho)$, where $0<m<n$, then
\begin{equation}\label{eq:32}
\frac{h_{j+1}}{k_{j+1}}=\frac{h_j+h_{j+2}}{\gcd(h_j+h_{j+2},k_j+k_{j+2})}
\biggm/\frac{k_j+k_{j+2}}{\gcd(h_j+h_{j+2},k_j+k_{j+2})}\ .
\end{equation}
\end{itemize}
\end{prop}

\begin{proof}
(i) There is nothing to prove if $m\in\{0,n\}$. If $0<m<n$ then,
in terms of Proposition~\ref{prop:6}(i), we have
$h_j=x_0+t^{\ast}h_{j+1}$, $k_j=y_0+t^{\ast}k_{j+1}$, and we
obtain $k_j h_{j+1}-k_{j+1}
h_j=(y_0+t^{\ast}k_{j+1})h_{j+1}-k_{j+1}(x_0+t^{\ast}h_{j+1})=y_0
h_{j+1}-x_0
k_{j+1}=\tfrac{x_0k_{j+1}+1}{h_{j+1}}h_{j+1}-x_0k_{j+1}=1$.

(ii) First, to see that $k_j h_{j+1}-h_jk_{j+1}=1$ and
$k_{j+1}h_{j+2}-h_{j+1}k_{j+2}=1$, apply assertion~(i) to each of
the pairs $\frac{h_j}{k_j}<\frac{h_{j+1}}{k_{j+1}}$ and
$\frac{h_{j+1}}{k_{j+1}}<\frac{h_{j+2}}{k_{j+2}}$. We have
$h_j=\tfrac{k_jh_{j+1}-1}{k_{j+1}}$,
$h_{j+2}=\tfrac{h_{j+1}k_{j+2}+1}{k_{j+1}}$, so then
$h_j+h_{j+2}=\tfrac{h_{j+1}}{k_{j+1}}(k_j+k_{j+2})$, and the
assertion follows.
\end{proof}

The following proposition is a tool of recurrent constructing
Farey subsequences~(\ref{eq:26}). In practice, such calculations
can be performed, for example, based on the successive fractions
mentioned in Remarks~{\rm\ref{rem:5}} and~\ref{rem:3}.

\begin{prop}
Let
$\frac{h_j}{k_j}<\frac{h_{j+1}}{k_{j+1}}<\frac{h_{j+2}}{k_{j+2}}$
be three successive fractions of
$\mathcal{F}(\mathbb{B}(n),m;\rho)$, where $0<m<n$.

\begin{itemize}
\item[\rm(i)] The integers $h_j$ and $k_j$ are computed by
\begin{align*}
h_j&=\left\lfloor\min\left\{\frac{h_{j+2}+m}{h_{j+1}},
\frac{k_{j+2}+n}{k_{j+1}},\frac{k_{j+2}
-h_{j+2}+n-m}{k_{j+1}-h_{j+1}}\right\}\right\rfloor
h_{j+1}-h_{j+2}\ ,\\
k_j&=\left\lfloor\min\left\{\frac{h_{j+2}+m}{h_{j+1}},
\frac{k_{j+2}+n}{k_{j+1}},\frac{k_{j+2}
-h_{j+2}+n-m}{k_{j+1}-h_{j+1}}\right\}\right\rfloor
k_{j+1}-k_{j+2}\ .
\end{align*}

\item[\rm(ii)] The integers $h_{j+2}$ and $k_{j+2}$ are computed by
\begin{align*}
h_{j+2}&=\left\lfloor\min\left\{\frac{h_j+m}{h_{j+1}},
\frac{k_j+n}{k_{j+1}},\frac{k_j
-h_j+n-m}{k_{j+1}-h_{j+1}}\right\}\right\rfloor h_{j+1}-h_j\ ,\\
k_{j+2}&=\left\lfloor\min\left\{\frac{h_j+m}{h_{j+1}},
\frac{k_j+n}{k_{j+1}},\frac{k_j
-h_j+n-m}{k_{j+1}-h_{j+1}}\right\}\right\rfloor k_{j+1}-k_j\ .
\end{align*}
\end{itemize}
\end{prop}

\begin{proof} To prove assertion~(i), note that, with respect to
Proposition~\ref{prop:7}(ii) and description~(\ref{eq:26}), we
have
\begin{equation*}
\begin{split}
\gcd(h_j+h_{j+2},k_j+k_{j+2})\cdot h_{j+1}&=h_j+h_{j+2}\leq
m+h_{j+2}\ ,\\ \gcd(h_j+h_{j+2},k_j+k_{j+2})\cdot
k_{j+1}&=k_j+k_{j+2}\leq n+k_{j+2}\ ,\\
\gcd(h_j+h_{j+2},k_j+k_{j+2})\cdot(k_{j+1}-h_{j+1})&=
(k_j-h_j)+(k_{j+2}-h_{j+2})\\&\phantom{=}\leq
(n-m)+(k_{j+2}-h_{j+2})\ ,
\end{split}
\end{equation*}
from where it follows that
\begin{multline*}
\gcd(h_j+h_{j+2},k_j+k_{j+2})\\=\left\lfloor\min\left\{\frac{h_{j+2}+m}{h_{j+1}},
\frac{k_{j+2}+n}{k_{j+1}},\frac{k_{j+2}
-h_{j+2}+n-m}{k_{j+1}-h_{j+1}}\right\}\right\rfloor\ ,
\end{multline*}
and we are done.

Assertion~(ii) is proved in an analogous way.
\end{proof}

{\small
\begin{rem}
For all elements $a\in\mathbb{V}_q(n)$ of rank $m:=\rho(a)$, the
Farey subsequences $\mathcal{F}(\mathbb{V}_q(n),a;\rho)$ are the
same, and we write $\mathcal{F}(\mathbb{V}_q(n),$ $m;\rho)$
instead of $\mathcal{F}(\mathbb{V}_q(n),a;\rho)$.

We have $\mathcal{F}(\mathbb{V}_q(n),m;\rho)=
\mathcal{F}(\mathbb{B}(n),m;\rho)$, for all $m$, $0\leq m\leq n$.
See also Remark~{\rm\ref{rem:2}}.
\end{rem}
}

{\small
\begin{rem}\label{rem:6}
In the present paper, we do not deal with the ascending Farey subsequences of the form
$\left(\tfrac{h}{k}\in\mathcal{F}_n:\ h\leq m\right)$, where $0< m\leq n$ {\rm(}see~{\rm Acketa and
\v{Z}uni\'{c}, 1991}, and Example~{\rm\ref{exmpl:1}}{\rm)}, including the classical Farey sequences
$\mathcal{F}_n$; see Section~{\rm\ref{SectionOnFareySubsequences}} for some references on $\mathcal{F}_n$.
Nevertheless, such Farey subsequences may be of use for the reader, and we list their basic properties.

\begin{itemize}
\item[\rm(i)]
In $\left(\tfrac{h}{k}\in\mathcal{F}_n:\ h\leq m\right)$, we have
\begin{multline*}
f_0=\tfrac{0}{1},\ f_1=\tfrac{1}{n},\\ f_{\left|\left(h/k\in\mathcal{F}_n:\ h\leq m\right)\right|-2}=
\tfrac{\min\{m,n-1\}}{\min\{m,n-1\}+1},\ f_{\left|\left(h/k\in\mathcal{F}_n:\ h\leq
m\right)\right|-1}=\tfrac{1}{1}.
\end{multline*}

\item[\rm(ii)]
\begin{itemize}
\item[\rm(a)]
If $f_t\in\left(\tfrac{h}{k}\in\mathcal{F}_n:\ h\leq
m\right)-\left\{\tfrac{1}{1}\right\}$ then
\begin{align*}
t&=\sum_{j\in\lbrack 1,n\rbrack}\phi \left( j;\bigl\lbrack
1,\min\{m,\lfloor j\cdot f_t\rfloor\}\bigr\rbrack \right)\\&=
-1+\sum_{d\geq 1}
\overline{\mu}(d)\cdot\left(\left\lfloor\frac{n}{d}\right\rfloor+\sum_{j\in
\left\lbrack 1,\left\lfloor
n/d\right\rfloor\right\rbrack}\min\left\{\left\lfloor\frac{m}{d}\right\rfloor,\left\lfloor
j\cdot f_t\right\rfloor\right\}\right)\ .
\end{align*}
\item[\rm(b)]
The cardinality of the sequence
$\left(\tfrac{h}{k}\in\mathcal{F}_n:\ h\leq m\right)$ equals
\begin{align*}
1+\sum_{j\in\lbrack 1,n\rbrack}\phi(j;\lbrack
1,\min\{m,j\}\rbrack)&=1+ \sum_{j\in\lbrack 1,m\rbrack}
\phi(j;\lbrack 1,j\rbrack )+\sum_{j\in\lbrack
m+1,n\rbrack}\phi(j;\lbrack 1,m\rbrack)\\ &= 1+\sum_{d\geq 1
}\overline{\mu}(d)\cdot\left(\left\lfloor\frac{n}{d}\right\rfloor-\frac{1}{2}
\left\lfloor\frac{m}{d}\right\rfloor\right)\cdot\left\lfloor\frac{m}{d}+1
\right\rfloor\ .
\end{align*}
\end{itemize}

\item[\rm(iii)]
Let $\frac{h}{k}\in\left(\tfrac{i}{j}\in\mathcal{F}_n:\ i\leq
m\right)$, $\frac{0}{1}<\frac{h}{k}<\frac{1}{1}$.

\begin{itemize}
\item[\rm(a)]
Let $x_0$ be the integer such that $kx_0\equiv -1\pmod{h}$ and
$m-h+1\leq x_0\leq m$. Define integers $y_0$ and $t^{\ast}$ by
$y_0:=\frac{kx_0+1}{h}$ and $t^{\ast}:=
\left\lfloor\min\{\frac{m-x_0}{h},\frac{n-y_0}{k}\}\right\rfloor$.
The fraction $\frac{x_0+t^{\ast} h}{y_0+t^{\ast} k}$ precedes the
fraction $\frac{h}{k}$ in $\left(\tfrac{h}{k}\in\mathcal{F}_n:\
h\leq m\right)$.

\item[\rm(b)]
Let $x_0$ be the integer such that $kx_0\equiv 1\pmod{h}$ and
$m-h+1\leq x_0\leq m$. Define integers $y_0$ and $t^{\ast}$ by
$y_0:=\frac{kx_0-1}{h}$ and $t^{\ast}:=
\left\lfloor\min\{\frac{m-x_0}{h},\frac{n-y_0}{k}\}\right\rfloor$.
The fraction $\frac{x_0+t^\ast h}{y_0+t^{\ast} k}$ succeeds the
fraction $\frac{h}{k}$ in $\left(\tfrac{h}{k}\in\mathcal{F}_n:\
h\leq m\right)$.
\end{itemize}

\item[\rm(iv)]
\begin{itemize}
\item[\rm(a)]
If $\frac{h_j}{k_j}<\frac{h_{j+1}}{k_{j+1}}$ are two successive
fractions of $\left(\tfrac{h}{k}\in\mathcal{F}_n:\ h\leq m\right)$
then {\rm(\ref{eq:31})} holds.

\item[\rm(b)]
If
$\frac{h_j}{k_j}<\frac{h_{j+1}}{k_{j+1}}<\frac{h_{j+2}}{k_{j+2}}$
are three successive fractions of
$\bigl(\tfrac{h}{k}\in\mathcal{F}_n:\ h\leq m\bigr)$
then~{\rm(\ref{eq:32})} holds.

\item[\rm(c)] If $\frac{h_j}{k_j}<\frac{h_{j+1}}{k_{j+1}}<\frac{h_{j+2}}{k_{j+2}}$ are three
successive fractions of $\bigl(\tfrac{h}{k}\in\mathcal{F}_n:\
h\leq m\bigr)$ then the integers $h_j$, $k_j$, $h_{j+2}$ and
$k_{j+2}$ are computed in the following way:
\begin{align*}
h_j&=\left\lfloor\min\left\{\tfrac{h_{j+2}+m}{h_{j+1}},
\tfrac{k_{j+2}+n}{k_{j+1}}\right\}\right\rfloor h_{j+1}-h_{j+2}\
,\\ k_j&=\left\lfloor\min\left\{\tfrac{h_{j+2}+m}{h_{j+1}},
\tfrac{k_{j+2}+n}{k_{j+1}}\right\}\right\rfloor k_{j+1}-k_{j+2}\
,\\ h_{j+2}&=\left\lfloor\min\left\{\tfrac{h_j+m}{h_{j+1}},
\tfrac{k_j+n}{k_{j+1}}\right\}\right\rfloor h_{j+1}-h_j\ ,\\
k_{j+2}&=\left\lfloor\min\left\{\tfrac{h_j+m}{h_{j+1}},
\tfrac{k_j+n}{k_{j+1}}\right\}\right\rfloor k_{j+1}-k_j\ .
\end{align*}
\end{itemize}

\item[\rm(v)]
If $\tfrac{1}{k}\in\left(\tfrac{i}{j}\in\mathcal{F}_n:\ i\leq
m\right)$, where $n>1$, for some $k>1$, then the fraction
$\tfrac{m+\min\left\{0,\left\lfloor\tfrac{n-km-1}{k}\right\rfloor\right\}}
{k\cdot\left(m+\min\left\{0,\left\lfloor\tfrac{n-km-1}{k}\right\rfloor\right\}\right)+1}$
precedes $\tfrac{1}{k}$, and the fraction
$\tfrac{m+\min\left\{0,\left\lfloor\tfrac{n-km+1}{k}\right\rfloor\right\}}
{k\cdot\left(m+\min\left\{0,\left\lfloor\tfrac{n-km+1}{k}\right\rfloor\right\}\right)-1}$
succeeds $\tfrac{1}{k}$ in $\left(\tfrac{i}{j}\in\mathcal{F}_n:\
i\leq m\right)$.
\end{itemize}
\end{rem}
}

\section{Relatively $r$-blocking elements in graded posets}
\label{SectionOnApplicationOfFareySubsequences}

Let $\{a\}$ be a nontrivial one-element antichain in $P$. Given a
map $\omega$, define the map
\begin{equation*}
\mathfrak{f}_{P,a;\omega}:\ \{r\in\mathbb{Q}:\ 0\leq
r<1\}\to\mathcal{F}(P,a;\omega)
\end{equation*}
by
\begin{equation*}
r\mapsto\max\bigl\{f\in\mathcal{F}(P,a;\omega):\ f\leq r\bigr\}\ .
\end{equation*}

Given an element $a$ of $\mathbb{B}(n)$, of rank $m:=\rho(a)>0$,
we write $\mathfrak{f}_{\mathbb{B}(n),m;\rho}(r)$ instead of
$\mathfrak{f}_{\mathbb{B}(n),a;\rho}(r)$.

The following assertion follows immediately from
Proposition~\ref{prop:5}(i).
\begin{cor}
Let $m$ and $n$ be positive integers such that $m<n$. If
$\mathfrak{f}_{\mathbb{B}(n),m;\rho}(r)=f_t\in\mathcal{F}(\mathbb{B}(n),m;\rho)$,
then
\begin{equation*}
t= \sum_{j\in\lbrack 1,n\rbrack}\phi \left(
j;\bigl\lbrack\max\{1,j+\min\{m,\lfloor j\cdot r\rfloor\}-n\},
\min\{m,\lfloor j\cdot r\rfloor\}\bigr\rbrack \right)\ .
\end{equation*}
\end{cor}

The Farey subsequences $\mathcal{F}(P,a;\omega)$ are, in
particular, of use because, given a nontrivial antichain $A$ in
$P$ and a map $\omega$, we have
\begin{equation*}
\mathbf{I}_r(P,A;\omega)= \bigcap_{a\in A}\
\mathbf{I}_{\mathfrak{f}_{P,a;\omega}(r)}(P,a;\omega)\ ,
\end{equation*}
cf.~Proposition~\ref{prop:1}(i).

Given a fraction $f$, we denote by $\underline{f}$ the numerator
of $f$, and we denote by $\overline{f}$ its denominator.

Given a nontrivial antichain $A$ in $P$ and a map $\omega$, define
a set $\mathcal{D}_r(P,A;\omega)\subset\mathbb{P}$ in the
following way:
\begin{multline*}
\mathcal{D}_r(P,A;\omega):=\bigcap_{a\in A}\\
\biggl(\bigcup_{f\in\mathcal{F}(P,a;\omega):\
\mathfrak{f}_{P,a;\omega}(r)<f<\frac{1}{1}}
\Bigl\{s\cdot\overline{f}:\ 1\leq s\leq\min\Bigl\{
\Bigl\lfloor\omega(a)/\underline{f}\Bigr\rfloor,\
\Bigl\lfloor\omega(P)/\overline{f}\Bigr\rfloor\Bigr\}\Bigr\}\\ \cup\left\{\omega(e):\ e\in\mathfrak{I}(a)-\{\hat{0}\}\right\}\biggr)\ .
\end{multline*}

This set of positive integers allows us to give the following
comment to Proposition~\ref{prop:4}(i).

\begin{prop}\label{prop:8} Let $P$ and $\omega$ satisfy the condition:
for any elements $a',a''\in P$, it holds
\begin{equation*}
\omega\bigl(\{a'\}\wedge_{\vartriangle}\{a''\}\bigr)=\omega(a')\
\Longleftrightarrow\ a'\leq a''\ .
\end{equation*}
Let $A$ be a nontrivial antichain in $P$, and let $k\in\lbrack
1,\omega(P)\rbrack$. Suppose that $|\bigcap_{a\in A
}\mathfrak{I}(a)-\{\hat{0}_P\}|=0$. If
$k\not\in\mathcal{D}_r(P,A;\omega)$ then
$|\mathbf{I}_{r,k}(P,A;\omega)|=0$.
\end{prop}

{\small
\begin{exmpl} If $A$ is an antichain in
$\mathbb{B}(6)$ such that $\{\rho(a):\ a\in A\} =\{2,3\}$, then
$\mathcal{D}_{\frac{1}{2}}(\mathbb{B}(6),A;\rho)=\{1,2,3\}$. Thus, if
the set $\mathbf{I}_{\!\frac{1}{2}}(\mathbb{B}(6),A;\rho)$ is
nonempty and if\/
$\mathbf{I}_{\!\frac{1}{2}}(\mathbb{B}(6),A;\rho)\ni b$, then
either $\{b\}=\bigcap_{a\in
A}\mathfrak{I}(a)-\{\hat{0}_{\mathbb{B}(6)}\}$ and $b$ is of rank
one, or\/ $b$ is of rank three.
\end{exmpl}
}

The concluding statement of the paper is a refinement of
Proposition~\ref{prop:4}(ii). Recall that the numbers $\nu(\cdot)$
are defined by~(\ref{eq:20}).

\begin{thm}\label{thm:1}
Let $P$ be a graded poset. If $A$ is a nontrivial antichain in $P$
then, on the one hand,
\begin{equation}\label{eq:37}
\mathbf{I}_r(P,A;\rho)= \left(\bigcap_{a\in A
}\mathfrak{I}(a)-\{\hat{0}_P\}\right)\ \cup\
\bigcup_{k\in\mathcal{D}_r(P,A;\rho)}
\Bigl(P^{(k)}\cap\mathfrak{F}\bigl(\boldsymbol{\mathsf{b}}_{\nu(r\cdot
k)-1}(A)\bigr)\Bigr)\ .
\end{equation}
On the other hand,
\begin{multline}\label{eq:38}
\mathbf{I}_r(P,A;\rho) = \left(\bigcap_{a\in A
}\mathfrak{I}(a)-\{\hat{0}_P\}\right)\\ \cup\ \bigcap_{a\in A }\ \
\bigcup_{f\in\mathcal{F}(P,a;\rho):\
\mathfrak{f}_{P,a;\rho}(r)<f}\ \
\bigcup_{s\in\left\lbrack 1,\ \min\left\{
\left\lfloor\rho(a)/\underline{f}\right\rfloor,\ \left\lfloor
\rho(P)/\overline{f}\right\rfloor\right\}\right\rbrack:\
s\cdot\overline{f}\in\mathcal{D}_r(P,A;\rho)}\\ \biggl(\
P^{(s\cdot\overline{f})}\cap
\Bigl(\mathfrak{F}\bigl(\boldsymbol{\mathsf{b}}_{s\cdot\underline{f}-1}(a)\bigr)-
\mathfrak{F}\bigl(\boldsymbol{\mathsf{b}}_{s\cdot\underline{f}}(a)\bigr)\Bigr)\
\biggr)\ .
\end{multline}
\end{thm}

\begin{proof}
First, in both expressions~(\ref{eq:37}) and~(\ref{eq:38}) we
consider the component $\bigcap_{a\in A
}\mathfrak{I}(a)-\{\hat{0}_P\}$; if this component, that
corresponds to the terminal fraction~$\tfrac{1}{1}$ of the Farey
subsequences $\mathcal{F}(P,a;\rho)$, $a\in A$, is nonempty then
any element of the component is a relatively $r$-blocking element
for $A$ in $P$. Further, if $k\not\in\mathcal{D}_r(P,A;\rho)$ then
the set $\mathbf{I}_{r,k}(P,A;\rho)-\left(\bigcap_{a\in A
}\mathfrak{I}(a)-\{\hat{0}_P\}\right)$ is empty, see
Proposition~\ref{prop:8}. Equality~(\ref{eq:37}) now follows from
Proposition~\ref{prop:4}(ii).

Let $a\in A$. We have
\begin{multline*}
\mathbf{I}_r(P,a;\rho)= \bigcup_{f\in\mathcal{F}(P,a;\rho):\
\mathfrak{f}_{P,a;\rho}(r)<f}\ \ \bigcup_{1\leq s\leq\min\left\{
\left\lfloor\rho(a)/\underline{f}\right\rfloor,\
\left\lfloor\rho(P)/\overline{f}\right\rfloor \right\}}\\
\bigl\{b\in P:\ \rho(b)=s\cdot \overline{f},\
\rho\bigl(\{b\}\wedge_{\vartriangle}\{a\}\bigr)= s\cdot
\underline{f}\bigr\}\ .
\end{multline*}
Further, we by~(\ref{eq:19}) have
\begin{equation*}
\begin{split}
\left\{b\in P:\ \rho(b)=s\cdot
\overline{f}\right\}&=P^{(s\cdot\overline{f})}\ ,\\ \bigl\{b\in
P:\ \rho\bigl(\{b\}\wedge_{\vartriangle}\{a\}\bigr)= s\cdot
\underline{f}\bigr\}&=\mathfrak{F}\bigl(\boldsymbol{\mathsf{b}}_{s\cdot
\underline{f}-1}(a)\bigr)-\mathfrak{F}\bigl(\boldsymbol{\mathsf{b}}_{s\cdot
\underline{f}}(a)\bigr)\ ,
\end{split}
\end{equation*}
and~(\ref{eq:38}) follows, with respect to
Proposition~\ref{prop:1}(i).
\end{proof}

{\small
\section*{References}

C.M.~Ablow and D.J.~Kaylor, ``A committee solution of the pattern
recognition problem,'' IEEE Trans. Inform. Theory vol.~11, no.~3,
pp.~453--455, 1965.

C.M.~Ablow and D.J.~Kaylor, ``Inconsistent homogeneous linear
inequalities,'' Bull. Amer. Math. Soc. vol.~71, no.~5, p.~724,
1965.

D.~Acketa and J.~\v{Z}uni\'{c}, ``On the number of linear
partitions of the $(m,n)$-grid,'' Inform. Process. Lett. vol.~38,
no.~3, pp.~163--168, 1991.

M.~Aigner, Combinatorial Theory, Grundlehren der Mathematischen
Wissenschaften, vol.~234, Springer-Verlag: Berlin, 1979.

G.E.~Andrews, ``Applications of basic hypergeometric functions,''
SIAM Review vol.~16, no.~4, pp.~283--294, 1974.

L.J.~Billera and A.~Bj\"{o}rner, ``Face numbers of polytopes and
complexes,'' in Handbook of Discrete and Computational Geometry,
J.E.~Goodman and J.~O'Rourke (eds.) CRC Press: Boca Raton,
New~York, pp.~291--310, 1997.

A.~Bj\"{o}rner,\hfill ``Topological methods,''\hfill in\hfill Handbook\hfill of\hfill
Combinatorics, \\ R.L.~Graham, M.~Gr\"{o}tschel and L.~Lov\'{a}sz
(eds.) Vol.~2,  Elsevier: Amsterdam, pp.~1819--1872, 1995.

A.~Bj\"{o}rner, L.M.~Butler and A.O.~Matveev, ``Note on a
combinatorial application of Alexander duality,''
J.~Combin.~Theory~Ser.~A vol.~80, no.~1, pp.~163--165, 1997.

A.~Bj\"{o}rner and A.~Hultman, ``A note on blockers in posets,''
Ann.~Comb. vol.~8, pp.~123--131, 2004.

A.~Bj\"{o}rner, M.~Las~Vergnas, B.~Sturmfels, N.~White and
G.M.~Ziegler, Oriented Matroids, Encyclopedia of Mathematics,
vol.~46, Cambridge University Press: Cambridge, 1993. Second
edition 1999.

A.~Bj\"{o}rner, I.~Peeva and J.~Sidman, ``Subspace arrangements
defined by products of linear forms,'' J.~London~Math.~Soc.~(2)
vol.~71, pp.~273--288, 2005.

W.~Bruns and J.~Herzog, Cohen-Macaulay Rings, Second edition,
Cambridge Studies in Advanced Mathematics, vol.~39, Cambridge
University Press: Cambridge, 1998.

A.A.~Buchstab, Teoria Chisel (in Russian) [Number Theory]
Uchpedgiz: Moscow, 1960.

V.M.~Buchstaber and T.E.~Panov, Toricheskie Deistviya v Topologii
i Kombinatorike (in Russian) [Torus Actions in Topology and
Combinatorics] Moskovskii Tsentr Nepreryvnogo Matematicheskogo
Obrazovaniya: Moscow, 2004.

C.J.~Colbourn, J.H.~Dinitz and D.R.~Stinson, ``Quorum systems
constructed from combinatorial designs,'' Inform. and Comput.
vol.~169, no.~2, pp.~160--173, 2001.

R.~Cordovil, K.~Fukuda and M.L.~Moreira, ``Clutters and
matroids,'' Discrete~Math. vol.~89, no.~2, pp.~161--171, 1991.

G.~Cornu\'{e}jols, Combinatorial Optimization. Packing and
Covering, CBMS-NSF Regional Conference Series in Applied
Mathematics, vol.~74, SIAM: Philadelphia, PA, 2001.

Y.~Crama and P.L.~Hammer, with contributions by C.~Benzaken,
J.C.~Bioch, E.~Boros, N.~Brauner, M.C.~Golumbic, V.~Gurvich,
L.~Hellerstein, T.~Ibaraki, A.~Kogan, K.~Makino, B.~Simeone and
B.~Vettier, Boolean Functions, book in preparation.

R.O.~Duda, P.E.~Hart and D.G.~Stork, Pattern Classification.
Second edition, Wiley-Interscience: New~York, 2001.

Z.~F\"{u}redi, ``Matchings and covers in hypergraphs,''
Graphs~Combin. vol.~4, no.~2, pp.~115--206, 1988.

R.L.~Graham, D.E.~Knuth and O.~Patashnik, Concrete Mathematics. A
Foundation for Computer Science, Second edition, Addison-Wesley:
Reading Massachusetts, 1994.

C.~Greene, ``The M\"{o}bius function of a partially ordered set,''
in Ordered Sets, Banff, Alta., 1981, NATO Adv. Study Inst. Ser. C:
Math. Phys. Sci., vol.~83, Reidel: Dordrecht-Boston, Mass.,
pp.~555--581, 1982.

M.~Gr\"{o}tschel, L.~Lov\'{a}sz and  A.~Schrijver, Geometric
Algorithms and Combinatorial Optimization, Second edition,
Algorithms and Combinatorics, vol.~2 Springer-Verlag: Berlin,
1993.

G.H.~Hardy and E.M.~Wright, An Introduction to the Theory of
Numbers, Fifth edition, Clarendon Press: Oxford, 1979.

T.~Hibi, Algebraic Combinatorics on Convex Polytopes, Carslaw
Publications: Glebe, Australia, 1992.

M.Yu.~Khachai, Komitetnye resheniya nesovmestnykh sistem
ogranichenii i metody obucheniya raspoznavaniyu (in Russian)
[Committee solutions of infeasible systems of constraints and
learning theory], DSci Thesis, Institute of Mathematics and
Mechanics, Russian Academy of Sciences, Ural Division,
Ekaterinburg, 2004.

M.Yu.~Khachai, Vl.D.~Mazurov and A.I.~Rybin, ``Committee
constructions for solving problems of selection, diagnostics, and
prediction,'' Proc. Steklov~Inst.~Math. Suppl.~1, pp.~S67--S101,
2002.

J.C.~Lagarias and C.P.~Tresser, ``A walk along the branches of the
extended Farey tree,'' IBM~J.~Res.~Develop. vol.~39, no.~3,
pp.~283--294, 1995.

D.E.~Loeb and A.R.~Conway, ``Voting fairly: transitive maximal
intersecting families of sets. In memory of Gian-Carlo Rota,''
J.~Combin.~Theory~Ser.~A vol.~91, no.~1-2, pp.~386--410, 2000.

A.O.~Matveev, ``A note on operators of deletion and contraction
for antichains,'' Int.~J.~Math.~Math.~Sci. vol.~31, no.~12,
pp.~725--730, 2002.

A.O.~Matveev, ``Extended blocker, deletion, and contraction maps
on antichains,'' Int.~J.~Math.~Math.~Sci. vol.~2003, no.~10,
pp.~607--616, 2003.

A.O.~Matveev, ``On blockers in bounded posets,''
Int.~J.~Math.~Math.~Sci. vol.~26, no.~10, pp.~581--588, 2001.

Vl.D.~Mazurov, Metod Komitetov v Zadachakh Optimizatsii i
Klassifikatsii (in Russian) [The Committee Method in Optimization
and Classification Problems], Nauka: Moscow, 1990.

Vl.D.~Mazurov, V.S.~Kazantsev, N.G.~Beletskii, A.I.~Krivonogov and
A.I.~Smirnov, ``Voprosy obosnovaniya i primeneniya komitetnykh
algoritmov raspoznavaniya'' (in Russian) [``Questions of the
justification and application of committee pattern recognition
algorithms''] in Raspoznavanie. Kassifikatsiya. Prognoz [Pattern
Recognition. Classification. Prediction], Yu.I.~Zhuravlev (ed.)
vol.~1, Nauka: Moscow, pp.~114--148, 1989.

Vl.D.~Mazurov and M.Yu.~Khachai, ``Komitetnye konstruktsii'' (in
Russian) [``Committee constructions''], Izv. Ural. Gos. Univ. Mat.
Mekh. vol~2, no.~14, pp.~77--108, 1999.

Vl.D.~Mazurov and M.Yu.~Khachai, ``Committees of systems of linear
inequalities,'' Automation and Remote Control vol.~65, no.~2,
pp.~193--203, 2004; translated from ``Komitety sistem lineinykh
neravenstv'' (in Russian), Avtomatika i Telemekhanika no.~2,
pp.~43--54, 2004.

E.~Miller and B.~Sturmfels, Combinatorial Commutative Algebra,
Graduate Texts in Mathematics, vol.~227, Springer-Verlag:
New~York, 2004.

M.~Naor and A.~Wool, ``The load, capacity, and availability of
quorum systems,'' SIAM~J.~Comput. vol.~27, no.~2, pp.~423--447,
1998.

P.~Orlik and H.~Terao, Arrangements of Hyperplanes, Grundlehren
der Mathematischen Wissenschaften [Fundamental Principles of
Mathematical Sciences], vol.~300, Springer-Verlag: Berlin, 1992.

C.E.~P\v{a}tra\c{s}cu and M.~P\v{a}tra\c{s}cu, ``Computing order
statistics in the Farey sequence,'' in Proc. 6th Algorithmic
Number Theory Symposium (ANTS 2004), Burlington, VT, USA, June
13-18, 2004, D.~Buell (ed.), Lecture Notes in Computer Science,
vol.~3076, Springer-Verlag: Heidelberg, pp.~358--366, 2004.

D.A.~Smith, ``Incidence functions as generalized arithmetic
functions.  I, II, III,'' Duke~Math.~J. vol.~34, pp.~617--633,
1967; vol.~36, pp.~15--30, 1969; vol.~36, pp.~353--367, 1969.

D.A.~Smith, ``Multiplication operators on incidence algebras,''
Indiana Univ.~Math.~J. vol.~20, pp.~369--383, 1970/1971.

R.P.~Stanley, Combinatorics and Commutative Algebra, Second
edition, Progress in Mathematics, Vol.~41, Birkhauser Boston,
Inc.: Boston, MA, 1996.

R.P.~Stanley, Enumerative Combinatorics, Vol. 1, Second edition,
Cambridge Studies in Advanced Mathematics, vol.~49, Cambridge
University Press: Cambridge, 1997.

T.~Zaslavsky, Facing up to Arrangements: Face-Count Formulas for
Partitions of Space by Hyperplanes, Mem. Amer. Math. Soc.,
no.~154, 1975.

G.M.~Ziegler, Lectures on Polytopes, Second edition, Graduate
Texts in Mathematics, vol.~152, Springer-Verlag: New~York, 1998.
}

\newpage

\section*{{\bf CORRIGENDUM:}\\ $\quad$\\
A.O.~Matveev, {\em Relative Blocking in Posets}, \\ Journal of Combinatorial Optimization, \\ {\bf 13} (2007), no.~4, 379--403\\ $\quad$}
\thispagestyle{empty}

\vspace{1cm}

{\em The use of the set of integers $\mathcal{D}_r(P,A;\omega)$ whose definition is given at the top of page~{\rm 401} may lead to situations where, mistakenly, some relatively $r$-blocking elements from the subposets $\mathbf{I}_r(P,A;\omega)\cap\mathfrak{I}(A)$ will not be taken into account. We correct the following inaccuracies}:

\begin{itemize}
\item
The definition of the set $\mathcal{D}_r(P,A;\omega)$ at the top of page~401
should read:
\begin{multline*}
\mathcal{D}_r(P,A;\omega):=\bigcap_{a\in A}\\
\biggl(\bigcup_{f\in\mathcal{F}(P,a;\omega):\
\mathfrak{f}_{P,a;\omega}(r)<f<\frac{1}{1}}
\Bigl\{s\cdot\overline{f}:\ 1\leq s\leq\min\Bigl\{
\Bigl\lfloor\omega(a)/\underline{f}\Bigr\rfloor,\
\Bigl\lfloor\omega(P)/\overline{f}\Bigr\rfloor\Bigr\}\Bigr\}\\ \cup\left\{\omega(e):\ e\in\mathfrak{I}(a)-\{\hat{0}\}\right\}\biggr)\ .
\end{multline*}

\item
The description of the set $\mathcal{D}_{\frac{1}{2}}(P,A;\rho)$ in Example~8.3 on page~401 should read: $\mathcal{D}_{\frac{1}{2}}(P,A;\rho)=\{1,2,3\}$.

\item
Expression~(8.2) of Theorem~8.4 on page~401 should read:

\begin{multline*}
\mathbf{I}_r(P,A;\rho) = \left(\bigcap_{a\in A
}\mathfrak{I}(a)-\{\hat{0}_P\}\right)\\ \cup\ \bigcap_{a\in A }\ \
\bigcup_{f\in\mathcal{F}(P,a;\rho):\
\mathfrak{f}_{P,a;\rho}(r)<f}\ \
\bigcup_{s\in\left\lbrack 1,\ \min\left\{
\left\lfloor\rho(a)/\underline{f}\right\rfloor,\ \left\lfloor
\rho(P)/\overline{f}\right\rfloor\right\}\right\rbrack:\
s\cdot\overline{f}\in\mathcal{D}_r(P,A;\rho)}\\ \biggl(\
P^{(s\cdot\overline{f})}\cap
\Bigl(\mathfrak{F}\bigl(\boldsymbol{\mathsf{b}}_{s\cdot\underline{f}-1}(a)\bigr)-
\mathfrak{F}\bigl(\boldsymbol{\mathsf{b}}_{s\cdot\underline{f}}(a)\bigr)\Bigr)\
\biggr)\ .\ \ \ \ \ \ \ \ \ (8.2)
\end{multline*}

\end{itemize}

\newpage

\end{document}